\documentclass[12pt,a4paper,twoside]{article}

\newcommand{\pageformat}[6]{\setlength{\hoffset}{-1in}
                  \setlength{\voffset}{-1in}
                  \addtolength{\hoffset}{#5}
                            \addtolength{\voffset}{#6}
                            \setlength{\oddsidemargin}{#1}
                            \setlength{\evensidemargin}{#2}
                            \setlength{\textwidth}{\paperwidth}
                  \addtolength{\textwidth}{-\oddsidemargin}
                  \addtolength{\textwidth}{-\evensidemargin}
                  \addtolength{\textwidth}{-\marginparsep}
                  \addtolength{\textwidth}{-\marginparwidth}
                            \setlength{\topmargin}{#3}
                            \setlength{\textheight}{\paperheight}
                  \addtolength{\textheight}{-\topmargin}
                  \addtolength{\textheight}{-\headheight}
                  \addtolength{\textheight}{-\headsep}
                  \addtolength{\textheight}{-\footskip}
                  \addtolength{\textheight}{-#4}}
\pageformat{2cm}{3cm}{25mm}{25mm}{1pt}{0pt}

\usepackage{ifthen}
\newboolean{article}
    \setboolean{article}{true}
\newboolean{report}
\newboolean{book}
\newboolean{letter}
\newboolean{german}
\newboolean{italian}
\newboolean{nobaselinestretch}
\newboolean{nosectionappendix}
\newboolean{oldtoc}
\newboolean{nosectionequn}
\newboolean{notheorem}

\ifthenelse{\boolean{german}}{
    \usepackage{german}}{}

\usepackage[latin1]{inputenc}

\usepackage{amsmath}
\usepackage{amssymb}
\usepackage[mathscr]{eucal}

\ifthenelse{\boolean{notheorem}}{}{
    \usepackage{theorem}}

\ifthenelse{\boolean{nobaselinestretch}}{}{
    \renewcommand{\baselinestretch}{1.25}}

\newenvironment{env}[2]{\begin{#1}#2\end{#1}}{}
    \newcommand{\beq}[1]{\begin{env}{equation}{#1}}
    \newcommand{\beqn}[1]{\begin{env}{equation*}{#1}}
    \newcommand{\bal}[1]{\begin{env}{align}{#1}}
    \newcommand{\baln}[1]{\begin{env}{align*}{#1}}
    \newcommand{\bga}[1]{\begin{env}{gather}{#1}}
    \newcommand{\bgan}[1]{\begin{env}{gather*}{#1}}
    \newcommand{\bflal}[1]{\begin{env}{flalign}{#1}}
    \newcommand{\bflaln}[1]{\begin{env}{flalign*}{#1}}
    \newcommand{\bmu}[1]{\begin{env}{multline}{#1}}
    \newcommand{\bmun}[1]{\begin{env}{multline*}{#1}}
    \newcommand{\bsp}[1]{\begin{env}{split}{#1}}

    \newcommand{\eeq}{\end{env}}
    \newcommand{\eeqn}{\end{env}}
    \newcommand{\eal}{\end{env}}
    \newcommand{\ealn}{\end{env}}
    \newcommand{\ega}{\end{env}}
    \newcommand{\egan}{\end{env}}
    \newcommand{\eflal}{\end{env}}
    \newcommand{\eflaln}{\end{env}}
    \newcommand{\emu}{\end{env}}
    \newcommand{\emun}{\end{env}}
    \newcommand{\esp}{\end{env}}

\newcommand{\lf}{\vspace{2ex}}

\renewcommand{\bf}[1]{\textbf{#1}}
\renewcommand{\it}[1]{\textit{#1}}
\renewcommand{\sc}[1]{\textsc{#1}}
\renewcommand{\sf}[1]{\textsf{#1}}

\renewcommand{\tt}[1]{\texttt{#1}}
\newcommand{\hl}[1]{\bf{\it{#1}}}

\newcommand{\mbf}[1]{\mathbf{#1}}
\newcommand{\msf}[1]{\text{\small$\sf{#1}$}}

\renewcommand{\mit}[1]{\mathit{#1}}
\newcommand{\cmc}[1]{\mathcal{#1}}
\newcommand{\eus}[1]{\mathscr{#1}}
\newcommand{\euf}[1]{\mathfrak{#1}}
\newcommand{\bb}[1]{\mathbb{#1}}

\newcommand{\mscriptsize}[1]{{\setlength{\arraycolsep}{.3ex}\text{\scriptsize$#1$}}}
\newcommand{\mtiny}[1]{{\setlength{\arraycolsep}{.3ex}\text{\tiny$#1$}}}
\newcommand{\nbd}[1]{$#1$\nobreakdash--}
\newcommand{\ol}[1]{\overline{#1}}

\newcommand{\vt}{\vartheta}

\newcommand{\vp}{\varphi}

\newcommand{\norm}[1]{\left\lVert#1\right\rVert}

\newcommand{\bfam}[1]{\bigl(#1\bigr)}
\newcommand{\Bfam}[1]{\Bigl(#1\Bigr)}
\newcommand{\AB}[1]{\langle#1\rangle}

\newcommand{\CB}[1]{\{#1\}}
\newcommand{\bCB}[1]{\bigl\{#1\bigr\}}
\newcommand{\BCB}[1]{\Bigl\{#1\Bigr\}}
\newcommand{\SB}[1]{[#1]}

\newcommand{\RO}[1]{[#1)}

\newcommand{\Matrix}[1]{\begin{pmatrix}#1\end{pmatrix}}

\newcommand{\sMatrix}[1]{\mscriptsize{\Matrix{#1}}}
\newcommand{\tMatrix}[1]{\mtiny{\Matrix{#1}}}

\newcommand{\sbar}[1]{\:\bar{#1}\:}
\newcommand{\sodot}{\sbar{\odot}}
\newcommand{\sbars}[1]{\:\bar{#1}^s\:}
\newcommand{\sodots}{\sbars{\odot}}
\newcommand{\set}[2][]{
    \ifthenelse{\equal{#1}{}}{
        \CB{#2}}{
        \CB{#1~|~#2}}}
\newcommand{\bset}[2][]{
    \ifthenelse{\equal{#1}{}}{
        \bCB{#2}}{
        \bCB{#1~|~#2}}}
\newcommand{\Bset}[2][]{
    \ifthenelse{\equal{#1}{}}{
        \BCB{#2}}{
        \BCB{#1~\big|~#2}}}

\DeclareMathOperator{\ls}{\normalfont\msf{span}}
\DeclareMathOperator{\cls}{\ol{\ls}}
\DeclareMathOperator*{\limind}{lim\,ind}

\DeclareMathOperator*{\soplus}{\ol{\bigoplus}^{\,\mit{s}}}

\DeclareMathOperator{\id}{\normalfont\msf{id}}

\renewcommand{\dim}{\operatorname{\msf{dim}}}
	\newcommand{\idim}{\operatorname{\msf{\scriptsize dim}}}

\newcommand{\C}{\bb{C}}

\newcommand{\E}{\bb{E}}

\newcommand{\N}{\bb{N}}

\newcommand{\R}{\bb{R}}
\newcommand{\bS}{\bb{S}}

\newcommand{\cA}{\cmc{A}}
\newcommand{\cB}{\cmc{B}}
\newcommand{\cC}{\cmc{C}}
\newcommand{\cD}{\cmc{D}}

\newcommand{\cO}{\cmc{O}}

\newcommand{\sB}{\eus{B}}

\newcommand{\sF}{\eus{F}}

\newcommand{\sK}{\eus{K}}

\newcommand{\sS}{\eus{S}}

\newcommand{\el}{\euf{l}}
\newcommand{\en}{\euf{n}}

\newcommand{\eH}{\euf{H}}

\newcommand{\U}{\mbf{1}}

\newcommand{\f}{\text{\scriptsize$\sF$}}
\newcommand{\s}{\text{\scriptsize$\sS$}}

\ifthenelse{\boolean{nosectionequn}}{}{
    \numberwithin{equation}{section}
    }

\ifthenelse{\boolean{article}\or\boolean{letter}\or\boolean{nosectionequn}}{
    \setboolean{nosectionappendix}{true}}{}
\ifthenelse{\boolean{nosectionappendix}}{}{
    \renewcommand{\appendix}{
        \chapter*{\appendixname}
        \addcontentsline{toc}{chapter}{\appendixname}
        \renewcommand{\thesection}{\Alph{section}}
        \setcounter{section}{0}}}
   
\ifthenelse{\boolean{report}\or\boolean{book}}{
    }{}

\ifthenelse{\boolean{notheorem}}{}{
        \newcommand{\mnname}{Mathematical note.}
        \newcommand{\enname}{End of the note.}
        \newcommand{\definame}{Definition.}
        \newcommand{\propname}{Proposition.}
        \newcommand{\lemname}{Lemma.}
        \newcommand{\exname}{Example.}
        \newcommand{\exername}{Exercise.}
        \newcommand{\remname}{Remark.}
        \newcommand{\obname}{Observation.}
        \newcommand{\thmname}{Theorem.}
        \newcommand{\corname}{Corollary.}
        \newcommand{\proofname}{Proof.}
    \ifthenelse{\boolean{german}}{
        \renewcommand{\mnname}{Mathematische Notiz.}
        \renewcommand{\enname}{Ende der Notiz.}
        \renewcommand{\exname}{Beispiel.}
        \renewcommand{\exername}{Übung.}
        \renewcommand{\remname}{Bemerkung.}
        \renewcommand{\obname}{Beobachtung.}
        \renewcommand{\thmname}{Satz.}
        \renewcommand{\corname}{Korollar.}
        \renewcommand{\proofname}{Beweis.}}{}
    \ifthenelse{\boolean{italian}}{
        \renewcommand{\mnname}{Nota matematica.}
        \renewcommand{\enname}{Fina della nota.}
        \renewcommand{\definame}{Definizione.}
        \renewcommand{\propname}{Proposizione.}
        \renewcommand{\exname}{Esempio.}
        \renewcommand{\exername}{Esercizio.}
        \renewcommand{\remname}{Nota.}
        \renewcommand{\obname}{Osservazione.}
        \renewcommand{\thmname}{Teorema.}
        \renewcommand{\corname}{Corollario.}
        \renewcommand{\proofname}{Dimostrazione.}

       \renewcommand{\appendixname}{Appendice}

       }{}
    \theoremheaderfont{\normalfont\bfseries}
    \theoremstyle{change}
        \theorembodyfont{\rmfamily}
            \newtheorem{emp}{}[section]
                \newcommand{\bemp}[1][]{
                    \begin{emp}\hskip-\labelsep\bf{#1}\hskip\labelsep}
                \newcommand{\eemp}{\end{emp}}
\newtheorem{itemp}[emp]{}
                \newcommand{\bitemp}[1][]{
                    \begin{itemp}\hskip-\labelsep\bf{#1}\hskip\labelsep\normalfont\itshape}
                \newcommand{\eitemp}{\end{itemp}}
            \newtheorem{mn}[emp]{\mnname}
                \newcommand{\bnm}{\begin{mn}~\begin{quotation}\renewcommand{\baselinestretch}{1}\small\noindent\ignorespaces}
                \newcommand{\enm}{\end{quotation}\hfill\bf{\enname}\end{mn}}
            \newtheorem{ex}[emp]{\exname}
                \newcommand{\bex}{\begin{ex}}
                \newcommand{\eex}{\end{ex}}
            \newtheorem{exer}[emp]{\exername}
                \newcommand{\bexer}{\begin{exer}}
                \newcommand{\eexer}{\end{exer}}
            \newtheorem{defi}[emp]{\definame}
                \newcommand{\bdefi}{\begin{defi}}
                \newcommand{\edefi}{\end{defi}}
            \newtheorem{rem}[emp]{\remname}
                \newcommand{\brem}{\begin{rem}}
                \newcommand{\erem}{\end{rem}}
            \newtheorem{ob}[emp]{\obname}
                \newcommand{\bob}{\begin{ob}}
                \newcommand{\eob}{\end{ob}}
        \theorembodyfont{\normalfont\itshape}
            \newtheorem{thm}[emp]{\thmname}
                \newcommand{\bthm}{\begin{thm}}
                \newcommand{\ethm}{\end{thm}}
            \newtheorem{prop}[emp]{\propname}
                \newcommand{\bprop}{\begin{prop}}
                \newcommand{\eprop}{\end{prop}}
            \newtheorem{cor}[emp]{\corname}
                \newcommand{\bcor}{\begin{cor}}
                \newcommand{\ecor}{\end{cor}}
            \newtheorem{lem}[emp]{\lemname}
                \newcommand{\blem}{\begin{lem}}
                \newcommand{\elem}{\end{lem}}
\newenvironment{empn}[1]{\lf\noindent\bf{#1}\ignorespaces\hskip\labelsep}{\lf}
		\newcommand{\bempn}[1]{\begin{empn}{#1}}
		\newcommand{\eempn}{\end{empn}}
		\newcommand{\bitempn}[1]{\begin{empn}{#1}\normalfont\itshape}
		\newcommand{\eitempn}{\end{empn}}
                \newcommand{\bnmn}{\begin{empn}{\mnname}~\begin{quotation}\renewcommand{\baselinestretch}{1}\small\noindent\ignorespaces}
                \newcommand{\enmn}{\end{quotation}\hfill\bf{\enname}\end{empn}}
		\newcommand{\bexn}{\begin{empn}{\exname}}
		\newcommand{\eexn}{\end{empn}}
		\newcommand{\bexern}{\begin{empn}{\exername}}
		\newcommand{\eexern}{\end{empn}}
		\newcommand{\bdefin}{\begin{empn}{\definame}}
		\newcommand{\edefin}{\end{empn}}
		\newcommand{\bremn}{\begin{empn}{\remname}}
		\newcommand{\eremn}{\end{empn}}
		\newcommand{\bobn}{\begin{empn}{\obname}}
		\newcommand{\eobn}{\end{empn}}

\newcommand{\qedsymbol}{~\rule[-0.35mm]{2mm}{2mm}}
    \newcounter{proof}[emp]
    \newenvironment{Proof}[1]{
        \vspace{1ex}
        \renewcommand{\item}[1][\stepcounter{proof}(\roman{proof})]%
            {##1\hskip\labelsep}
        \noindent\textsc{#1\hskip\labelsep}}{
        \nolinebreak\qedsymbol}
    \newcommand{\proof}[1][\proofname]{
        \begin{Proof}{#1}\ignorespaces}
    \newcommand{\qed}{\end{Proof}}
    \newcommand{\noqed}{
        \renewcommand{\qedsymbol}{}
        \end{Proof}}}
    \ifthenelse{\boolean{italian}}{
        \renewcommand{\proofname}{Dimostrazione.}}{}

\usepackage[varg]{txfonts}

\begin{document}

\title{Unit Vectors, Morita Equivalence and Endomorphisms}
\author{}
\author{
~\\
Michael Skeide\thanks{This work is supported by research fonds of the Department S.E.G.e S.\ of University of Molise.}\\\\
{\small\itshape Dipartimento S.E.G.e S.}\\
{\small\itshape Università degli Studi del Molise}\\
{\small\itshape Via de Sanctis}\\
{\small\itshape 86100 Campobasso, Italy}\\
{\small{\itshape E-mail: \tt{skeide@math.tu-cottbus.de}}}\\
{\small{\itshape Homepage: \tt{http://www.math.tu-cottbus.de/INSTITUT/lswas/\_skeide.html}}}\\\\
}
\date{December 2004, This revision September 2007}

{
\renewcommand{\baselinestretch}{1}
\maketitle


\vspace{14ex}
\begin{abstract}
\noindent
We solve two problems in the theory of correspondences that have important implications in the theory of product systems. The first problem is the question whether every correspondence is the correspondence associated (by the representation theory) with a unital endomorphism of the algebra of all adjointable operators on a Hilbert module. The second problem is the question whether every correspondence allows for a nondegenerate faithful representation on a Hilbert space. We also resolve an extension problem for representations of correspondences and we provide new efficient proofs of several well-known statements in the theory of representations of \nbd{W^*}algebras.
\end{abstract}

}




\section{Introduction}

{\parskip0.5ex plus 0.5ex minus 0.5ex
Let $\cB$ be a \nbd{C^*}algebra. With every unital strict endomorphism of the \nbd{C^*}algebra $\sB^a(F)$ of all adjointable operators on a Hilbert \nbd{\cB}module $F$ there is associated a correspondence $F_\vt$ over $\cB$ (that is, a Hilbert \nbd{\cB}bimodule) such that
\bal{\label{PScorr}
F
&~=~
F\odot F_\vt
&
\vt(a)
~=~
a\odot\id_{F_\vt}.
}\eal
In other words, $\vt$ is amplification of $\sB^a(F)$ with the \it{multiplicity correspondence} $F_\vt$. (This is just the representation theory of $\sB^a(F)$.) The same is true for a \nbd{W^*}module (where $\vt$  is normal and the tensor product is that of \nbd{W^*}correspondences.)
 
\lf\noindent
\bf{Problem 1.~}
Given a correspondence $E$ over a \nbd{C^*} (or \nbd{W^*})algebra $\cB$, construct a unital strict (or normal) endomorphism $\vt$ of some $\sB^a(F)$ such that $E$ is the \it{multiplicity correspondence} $F_\vt$ associated with $\vt$.

\lf
An intimately related problem (in the \nbd{W^*}case, in fact, an equivalent problem) is the following.

\lf\noindent
\bf{Problem 2.~}
Find a nondegenerate faithful (normal) representation of the (\nbd{W^*})correspondence $E$ over $\cB$ on some Hilbert space.

\lf
In these notes we resolve Problem 1 for strongly full \nbd{W^*}correspondences and for full correspondences over a unital \nbd{C^*}algebra. We resolve Problem 2 for correspondences and \nbd{W^*}correspondences that are \it{faithful} in the sense that the left action of the correspondence is faithful. (Recall that, by definition, all correspondences have nondegenerate left action.) The conditions, fullness for Problem 1 and faithfulness for Problem 2, are also necessary. So, apart from Problem 1 for a nonunital \nbd{C^*}algebra we present a complete solution of the two problems. We explain that in the \nbd{W^*}case the two problems are dual to each other in the sense of the \it{commutant} of von Neumann correspondences. Throughout, \it{en passant} we furnish a couple of new, simple proofs for known statements that illustrate how useful our methods are.

\lf
The study of representations of correspondences goes back, at least, to Pimsner \cite{Pim97} and, in particular, to Muhly and Solel \cite{MuSo98} and their forthcoming papers. Hirshberg \cite{Hir05a} resolved Problem 2 for \nbd{C^*}correspondences that are faithful and full. We add here (by furnishing a completely different proof) that the hypothesis of fullness is not necessary and that in the \nbd{W^*}case the representation can be chosen normal.

Problem 1 is the ``reverse'' of the representation theory of $\sB^a(F)$; Skeide \cite{Ske02,Ske03c,Ske05a} and Muhly, Skeide and Solel \cite{MSS06}.

Our interest in the solution of the Problems 1 and 2 has its common root in the theory of \nbd{E_0}semigroups (that is, semigroups of unital endomorphisms) of $\sB^a(F)$ and their relation with product systems of correspondences. Arveson \cite{Arv89} associated with every normal \nbd{E_0}semigroup on $\sB(H)$ ($H$ a Hilbert space) a product system of Hilbert spaces (\hl{Arveson system}, for short) that comes along with a natural faithful representation. Finding a faithful representation of a given Arveson system is equivalent to that this Arveson system is the one associated as in \cite{Arv89} with an \nbd{E_0}semigroup. In the three articles \cite{Arv90a,Arv89a,Arv90} Arveson showed that every Arveson system admits a faithful representation, that is, it is the Arveson system associated with an \nbd{E_0}semigroup as in \cite{Arv89}.

Bhat \cite{Bha96} constructed from a normal \nbd{E_0}semigroup on $\sB(H)$ a second Arveson system (the \hl{Bhat system} of the \nbd{E_0}semigroup) that turns out to be anti-isomorphic to the one constructed by Arveson \cite{Arv89}. The Bhat system is related to the endomorphisms of the \nbd{E_0}semigroup via Equation \eqref{PScorr}.

It is Bhat's point of view that generalizes directly to \nbd{E_0}semigroups of $\sB^a(F)$, while Arveson's point of view works only when $F$ is a von Neumann module. (In fact, the two product systems are no longer just anti-isomorphic, but as explained in Skeide \cite{Ske03c} they turn out to be \it{commutants} of each other; see Section \ref{comndrsec}.)

In Skeide \cite{Ske06} we presented a short and elementary proof of Arveson's result that every Arveson system is the one associated with an \nbd{E_0}semigroup. This proof uses essentially the fact that it is easy to resolve the problem for discrete time $t\in\N_0$ or, what is the same, for a single Hilbert space $\eH$ (that generates a discrete product system $\bfam{\eH^{\otimes n}}_{n\in\N_0}$). If we want to apply the idea of the proof in \cite{Ske06} also to Hilbert and von Neumann modules, then we must first resolve the problem for a single correspondence $E$ (that generates a discrete product system $\bfam{E^{\odot n}}_{n\in\N_0}$). This is precisely what we do in these notes: Solving Problem 1 means that $\bfam{E^{\odot n}}_{n\in\N_0}$ is the product system of the discrete \nbd{E_0}semigroup $\bfam{\vt^n}_{n\in\N_0}$. Solving Problem 2 means finding a faithful representation of the whole discrete product system $\bfam{E^{\odot n}}_{n\in\N_0}$. In fact, in the meantime we did already use the results of these notes (or ideas leading to them) to resolve the continuous time case for Hilbert modules \cite{Ske07,Ske06p6} and von Neumann modules \cite{Ske07p} (in preparation).

In the solution of Problems 1 and 2 the concepts of \it{unit vectors} in Hilbert or \nbd{W^*}modules and of \it{Morita equivalence} for (\nbd{W^*})correspondences and modules play a crucial role. In fact, if a correspondence $E$ has unit vector $\xi$ (that is, $\AB{\xi,\xi}=\U\in\cB$ so that, in particular, $E$ is full and $\cB$ is unital), then it is easy to construct a unital \it{endomorphism} $\vt$ on some $\sB^a(F)$ that has $E$ as associated multiplicity correspondence $F_\vt$; see Section \ref{psunisec}. Morita equivalence helps to reduce Problem 1 for (strongly) full (\nbd{W^*})correspondences to the case when $E$ has a unit vector. In fact, even if a (strongly) full $E$ does not have a unit vector, then \it{cum grano salis} (that is, up to suitable completion) the space of \nbd{E}valued matrices $M_\en(E)$ of sufficiently big dimension will have a unit vector. The correspondences $M_\en(E)$ and $E$ are Morita equivalent in a suitable sense, and in Theorem \ref{psMethm} we show that solving Problem 1 for $M_\en(E)$ is equivalent to solving Problem 1 for $E$ itself. Last but not least, we mention that Morita equivalence is at the heart of the representation theory of $\sB^a(F)$ which we use to determine the correspondence of an endomorphism; see Example \ref{repthex}.

Problem 2, instead, in the \nbd{W^*}case (Theorem \ref{W*Hirthm}) is a simple consequence of the well-known fact that two faithful normal nondegenerate representations of a \nbd{W^*}algebra have unitarily equivalent amplifications. In order to illustrate how simply this result can be derived making appropriate use of unit vectors and quasi orthonormal bases in von Neumann modules, we include a proof (Corollary \ref{ampisocor}). The \nbd{C^*}case (Theorem \ref{C*Hirthm}) is a slightly tedious reduction to the \nbd{W^*}case. In Theorem \ref{EE'switchthm} we show that the \nbd{W^*}versions of Problem 1 and Problem 2 are, actually, equivalent. However, while the \nbd{C^*}version of Problem 2 can be reduced to the \nbd{W^*}version, a similar procedure is not possible for Problem 1. (Given a full correspondence over a possibly nonunital \nbd{C^*}algebra $\cB$, we can resolve Problem 1 for the enveloping \nbd{W^*}correspondence over $\cB^{**}$. But, we do not know a solution to the problem how find a (strongly dense) \nbd{\cB}submodule $F$ of the resulting \nbd{\cB^{**}}module $F^{**}$ such that the endomorphism $\vt$ of $\sB^a(F^{**})$ restricts suitably to an endomorphism of $\sB^a(F)$.)

\lf
These notes are organized as follows. In Section \ref{psunisec} we explain the relation between \nbd{E_0}semi\-groups on $\sB^a(E)$ and product systems. We discuss a case in which it is easy to construct for a product system an \nbd{E_0}semigroup with which the product system is associated. In Observation \ref{P1ob} we explain how this leads to a simple solution of Problem 1 in the case when the correspondence has a unit vector.

In Section \ref{C*uvsec} we show that a finite multiple of a full Hilbert module over a unital \nbd{C^*}algebra has a unit vector (Lemma \ref{univeclem}). Apart from a simple consequence about finitely generated Hilbert modules (Corollary \ref{Kunicor}), this lemma is crucial for the solution of the \nbd{C^*}version of Problem 1 in Section \ref{C*sec}. In Section \ref{W*uvsec} we prove the \nbd{W^*}analogue of Lemma \ref{univeclem}, Lemma \ref{Wuniveclem}: A suitable multiple of a strongly full \nbd{W^*}module has a unit vector. The proof is considerably different from that of Lemma \ref{univeclem}. It makes use of quasi orthonormal bases. We use the occasion to illustrate how easily some basic facts about representations of von Neumann algebras, like the \it{amplification-induction theorem}, may be derived. Utilizing in an essential way Lemma \ref{Wuniveclem}, we give a simple proof of the well-known fact that faithful normal representations of a \nbd{W^*}algebra have unitarily equivalent amplifications (Corollary \ref{ampisocor}). A proof of that result is also included to underline how simple a self-contained proof of the solution to Problem 2 (Theorems \ref{W*Hirthm} and \ref{C*Hirthm}) actually is.

Section \ref{MoritaSEC} introduces the necessary notions of Morita equivalence. Apart from (strong) Morita equivalence for $C^*$--{~} and \nbd{W^*}algebras, we discuss Morita equivalence for correspondences (Muhly and Solel \cite{MuSo00}) and Morita equivalence for Hilbert and \nbd{W^*}modules (new in these notes). We state the obvious generalization of Morita equivalence for correspondences to product systems. Two full Hilbert modules have strictly isomorphic operator algebras, if and only if they are Morita equivalent. In that case, two endomorphisms (\nbd{E_0}semigroups on the isomorphic operator algebras are (cocycle) conjugate, if and only if they have Morita equivalent correspondences (product systems); see Proposition \ref{conprop} and Corollary \ref{concor}. The central result is Theorem \ref{psMethm}, which asserts that in the \nbd{W^*}case solvability of Problem 1 does not change under Morita equivalence.

In Sections \ref{mthmsec} and \ref{C*sec} we resolve Problem 1 for \nbd{W^*}correspondences (Theorem \ref{normunithm}) and correspondences over unital \nbd{C^*}algebras (Theorem \ref{strunithm}), respectively. While the \nbd{W^*}case runs smoothly after the preparation in Sections \ref{W*uvsec} and \ref{MoritaSEC}, in the \nbd{C^*}case we have to work considerably. In both sections we spend some time to explain where the difficulties in the \nbd{C^*}case actually lie.

Section \ref{repsec} contains the complete solution to Problem 2. Taking into account Corollary \ref{ampisocor}, the treatment is a self-contained. A simple consequence of Sections \ref{psunisec} and \ref{repsec} are Theorem \ref{AKthm} and its corollary, which assert that a faithful endomorphism is a restriction to a subalgebra of some inner endomorphism on $\sB(H)$. In Theorem \ref{ndthm} we solve the apparently open problem to find a nondegenerate \it{extension} to a normal faithful representation (in the language of Muhly and Solel \cite{MuSo98}, a fully coisometric extension of a \nbd{\sigma}continuous faithful isometric covariant representation) of a \nbd{W^*}correspondence.

In Section \ref{comndrsec} we show that the \nbd{W^*}version of Problem 1 and Problem 2 are equivalent under the \it{commutant} of von Neumann correspondences (Theorems \ref{EE'switchthm} and \ref{PSswitchthm}). The fact that, to that goal, we have to discuss the basics about von Neumann modules and von Neumann correspondences has the advantage that we provide also simple proofs for many statements about \nbd{W^*}modules, used earlier in these notes. As some more consequences of Corollary \ref{ampisocor} and the language used in Section \ref{comndrsec}, we furnish new proofs for the well-known results Corollary \ref{freecor} (a sort of \it{Kasparov absorption theorem} for \nbd{W^*}modules) and Corollary \ref{Mecor} (a couple of criteria for when two \nbd{W^*}algebras are Morita equivalent). Corollary \ref{freecor} is also the deeper reason for that the solutions to our Problems 1 and 2 in the \nbd{W^*}case may be chosen of a particularly simple form; see Observations \ref{uniexob} and \ref{repampob}.

In Section \ref{exsec} we discuss our results in two examples.

\lf\noindent
\bf{A note on the first version.~}
These notes are a very far reaching revision of the version of the preprint published as \cite{Ske04p}. The main results (Theorems \ref{normunithm}, \ref{strunithm}, and \ref{W*Hirthm}) and essential tools (Lemmata \ref{univeclem} and \ref{Wuniveclem}, Theorems \ref{psMethm} and \ref{EE'switchthm}) have been present already in \cite{Ske04p}. But while Theorem \ref{W*Hirthm} in \cite{Ske04p} has been proved by reducing it to Theorem \ref{normunithm} via the \it{commutant}, the new simple proof we give here is now independent of Section \ref{comndrsec} and Theorem \ref{normunithm}. New in this revision are the proof of Hirshberg's result \cite{Hir05a} that works also in the nonfull case (Theorem \ref{C*Hirthm}), and the extension result Theorem \ref{ndthm}. A couple of very simple proofs of well-known results has been included. Finally, the discussion of the examples in Section \ref{exsec} has been shortened drastically. For some details in these examples we find it convenient to refer the reader to the old version \cite{Ske04p}.

\lf\noindent
\bf{Acknowledgements.~}
We wish to express our gratitude to P.\ Muhly for warm hospitality during several stays in Iowa City with many inspiring discussions that influenced this work. We thank I.\ Hirshberg for intriguing discussions around his article \cite{Hir05a} and A.\ Skalski for making us aware of the extension problem solved in Theorem \ref{ndthm}. We are particularly grateful to B.\ Solel for pointing out to us an inaccuracy in the first version. Last but not least, we acknowledge also the support by Research Funds of the Department S.E.G.e S.\ of University of Molise.

\lf\noindent
\bf{Notations, conventions and some basic properties.~}

\bemp\label{Iop}
By $\sB^a(E)$ we denote the algebra of adjointable operators on a Hilbert \nbd{\cB}module $E$. A linear map $\vt\colon\sB^a(E)\rightarrow\sB^a(F)$ is \hl{strict}, if it is continuous on bounded subsets for the strict topologies of $\sB^a(E)$ and $\sB^a(F)$. Recall that a unital endomorphism $\vt$ of $\sB^a(E)$ is strict, if and only if the action of the compact operators $\sK(E)$ is already nondegenerate: $\cls\sK(E)E=E$. The \nbd{C^*}algebra of \hl{compact operators} is the completion $\sK(E):=\ol{\sF(E)}$ of the finite-rank operators, and the pre-\nbd{C^*}algebra of \hl{finite-rank operators} is the linear span $\sF(E):=\ls\bCB{xy^*\colon x,y\in E}$ of the \hl{rank-one operators} $xy^*\colon z\mapsto x\AB{y,z}$.
\eemp

\bemp\label{Ifulluni}
The \hl{range ideal} of a Hilbert \nbd{\cB}module is the closed ideal $\cB_E:=\cls\AB{E,E}$ in $\cB$. A Hilbert \nbd{\cB}module $E$ is \hl{full}, if $\cB_E=\cB$. A \hl{unit vector} in a Hilbert \nbd{\cB}module $E$ is an element $\xi\in E$ fulfilling $\AB{\xi,\xi}=\U\in\cB$. This means, in particular, that $\cB$ is unital and that $E$ is full.
\eemp

\bemp\label{Icorr}
A \hl{correspondence from $\cA$ to $\cB$} is a Hilbert \nbd{\cB}module with a nondegenerate\hl{(!)} left action of $\cA$. When $\cA=\cB$, we shall also say correspondence \hl{over $\cB$}. We say a correspondence from $\cA$ to $\cB$ is \hl{faithful}, if the left action of $\cA$ defines a faithful homomorphism. Every \nbd{C^*}algebra $\cB$ is a correspondence over itself, the \hl{trivial correspondence over $\cB$}, with inner product $\AB{b,b'}:=b^*b'$ and the natural bimodule operations. The \nbd{\cB}subcorrespondence of the trivial correspondence $\cB$ correspond precisely to the closed ideals.
\eemp

\bemp\label{IspecialE}
Every Hilbert \nbd{\cB}module is a correspondence from $\sB^a(E)$ to $\cB$ that may be viewed also as a correspondence from $\sK(E)$ (or any \nbd{C^*}algebra in between $\sK(E)$ and $\sB^a(E)$) to $\cB_E$ (or any \nbd{C^*}algebra in between $\cB_E$ and $\cB$). The \hl{dual correspondence} of $E$ is the correspondence $E^*=\bCB{x^*\colon x\in E}$ from $\cB$ to $\sB^a(E)$. It consists of mappings $x^*\colon y\mapsto\AB{x,y}$ in $\sB^a(E,\cB)$ with inner product $\AB{x^*,y^*}:=xy^*$ and bimodule operations $bx^*a:=(a^*xb^*)^*$. We note that $\sK(E^*)=\cB_E$ and that the range ideal is $\sB^a(E)_{E^*}=\sK(E)$. The left action of $\cB_E$ is, indeed, faithful so that $E^*$ may be viewed as faithful and full correspondence from $\cB_E$ to $\sK(E)$.
\eemp

\bemp\label{Itp}
The \hl{(internal) tensor product} of a correspondence $E$ from $\cA$ to $\cB$ and a correspondence $F$ from $\cB$ to $\cC$ is that unique correspondence $E\odot F$ from $\cA$ to $\cC$ that is generated by the range of a left \nbd{\cA}linear mapping $(x,y)\mapsto x\odot y$ fulfilling $\AB{x\odot y,x'\odot y'}=\AB{y,\AB{x,x'}y'}$.

For every correspondence $E$ from $\cA$ to $\cB$ we have the \hl{canonical identifications} $\cA\odot E\cong E$ via $a\odot x\mapsto ax$ (recall that, by our convention in Section \ref{Icorr}, $\cA$ acts nondegenerately), and $E\odot\cB\cong E$ via $x\odot b\mapsto xb$. One easily verifies that $x\odot y^*\mapsto xy^*$ defines an isomorphism $E\odot E^*\rightarrow\sK(E)$ of correspondences over $\sK(E)$ (or over $\sB^a(E)$). Similarly, $x^*\odot y\mapsto\AB{x,y}$ defines an isomorphism $E^*\odot E\rightarrow\cB_E$ of correspondences over $\cB_E$ (or over $\cB$). We will always identify these correspondences.
\eemp

\bemp\label{IW*}
A \hl{\nbd{W^*}module} is a Hilbert module over a \nbd{W^*}algebra that is self-dual. A Hilbert \nbd{\cB}module $E$ is \hl{self-dual}, if every bounded right-linear map $\Phi\colon E\rightarrow\cB$ has the form $x^*$ of some $x\in E$. Every Hilbert module over a \nbd{W^*}algebra admits a unique minimal self-dual extension \cite{Pas73,Rie74a,Ske00b,Ske05c}; see Remarks \ref{dscrem} and \ref{dctrem}. A \hl{\nbd{W^*}correspondence} $E$ is a \nbd{C^*}correspondence from a \nbd{W^*}algebra $\cA$ to a \nbd{W^*}algebra $\cB$ and a \nbd{W^*}module, such that all maps $\AB{x,\bullet x}\colon\cA\rightarrow\cB$ $(x\in E)$ are normal. The \hl{\nbd{W^*}tensor product} of a \nbd{W^*}correspondence $E$ from $\cA$ to $\cB$ and a \nbd{W^*}correspondence $F$ from $\cB$ to $\cC$ is the self-dual extension $E\sodots F$ of the tensor product $E\odot F$. This extension is a \nbd{W^*}correspondence from $\cA$ to $\cC$.

If $E$ is a \nbd{W^*}module over $\cB$, then the \hl{extended linking algebra} $\sMatrix{\cB&E^*\\E&\sB^a(E)}$ is a \nbd{W^*}algebra. By restriction, this equips every corner with a \hl{\nbd{\sigma}weak topology} and a \hl{\nbd{\sigma}strong topology}. Properties of these topologies in the linking algebra directly turn over to the corners. Consequently, we say a map $\eta\colon E\rightarrow F$ is \hl{normal}, if it is the restriction to the \nbd{2}\nbd{1}corners of a normal map between the extended linking algebras.

Most of the statements about \nbd{W^*}modules and \nbd{W^*}correspondences have considerably simpler proofs in the equivalent categories of von Neumann modules and von Neumann correspondences \cite{Rie74a,BDH88,Ske00b,Ske03c,Ske05c,Ske06b}; see Section \ref{comndrsec}.
\eemp

\bemp\label{Irep}
Suppose $E$ is a Hilbert \nbd{\cB}module. Let us choose a (nondegenerate) representation $\pi$ of $\cB$ on a Hilbert space $G$. We may construct the Hilbert space $H:=E\odot G$, and the \hl{induced representation} $\rho^\pi$ of $\sB^a(E)$ on $\sB(H)$ by setting $\rho^\pi(a):=a\odot\id_G$. We define the \hl{induced representation} $\eta^\pi\colon E\rightarrow\sB(G,H)$ of $E$ from $G$ to $H$ by setting $\eta^\pi(x)g=x\odot g$. (That is, $\eta^\pi$ fulfills $\eta^\pi(x)^*\eta^\pi(y)=\pi(\AB{x,y})$ and $\eta^\pi(xb)=\eta^\pi(x)\pi(b)$. Obviously, $\eta^\pi(ax)=\rho^\pi(a)\eta^\pi(x)$.) The maps $\pi$, $\eta^\pi$, $(\eta^\pi)^*:=*\circ\eta^\pi\circ*$, and $\rho^\pi$ give rise to the (nondegenerate) \hl{induced representation} $\Pi:=\sMatrix{\pi&(\eta^\pi)^*\\\eta^\pi&\rho^\pi}$ of the extended linking algebra on $G\oplus H$. So, all mappings are completely contractive.

In the language of von Neumann modules it is not difficult to show that for a (strongly full) \nbd{W^*}module the induced representation of the extended linking algebra is normal, if (and only if) $\pi$ is normal.

If $E$ is a correspondence from $\cA$ to $\cB$, then we will also speak of the \hl{induced representation $\rho^\pi_\cA\colon\cA\rightarrow\sB^a(E)\rightarrow\sB(H)$} of $\cA$ on $H$. If $E$ is faithful, then simply $\rho^\pi_\cA=\rho^\pi\upharpoonright\cA$.
\eemp

\bemp\label{Icard}
We will often need multiples of an arbitrary cardinality $\en$ of Hilbert spaces or modules. If $\en$ is a cardinal number, then we always assume that we have fixed a set $S$ with cardinality $\#S=\en$ so that $E^\en:=\bigoplus_{s\in S}E=\ol{\varprod_{s\in S}E}$ has a well-specified meaning. If $T$ is another set having that cardinality, then, by definition of cardinality, there exists a bijection between $S$ and $T$ that induces a \it{canonical} isomorphism from $\bigoplus_{s\in S}E$ to $\bigoplus_{t\in T}E$. So, $\C^\en$ is \it{the} Hilbert space (up to canonical isomorphism) of dimension $\en$. For every Hilbert space $\eH$ we may write $\eH=\C^{\idim\eH}$. Of course, $E^\en=E\otimes\C^\en$ (or $=\C^\en\otimes E$) in the sense of \it{external} tensor products, and we may write $E\otimes\eH=E^{\idim\eH}$. Amplifications $a\otimes\id_\eH$ of a map $a$ on $E$ in the tensor product picture, will be written as $a^{\idim\eH}$ in the direct sum picture.
\eemp
}


\section{Prerequisites on $E_0$--semigroups and product systems}\label{psunisec}

Let $\bS$ denote either the semigroup of nonnegative integers $\N_0=\CB{0,1,\ldots}$ or the semigroup of nonnegative reals $\R_+=\RO{0,\infty}$. In this section we explain the relation between a strict \nbd{E_0}semigroup $\vt=\bfam{\vt_t}_{t\in\bS}$ and its product system $E^\odot=\bfam{E_t}_{t\in\bS}$. In these notes we are mainly interested in the \hl{discrete} case $\bS=\N_0$. However, there is no reason to restrict the present discussion to the discrete case. In fact, many results we prove in these notes hold in the general case. They find their applications in Skeide \cite{Ske07,Ske06p6,Ske07p}, where we discuss several variants of the \hl{continuous time} case $\bS=\R_+$, and, in a different context, in Skeide \cite{Ske06p3}.

Let $E$ be a Hilbert module over a \nbd{C^*}algebra $\cB$ and let $\vt=\bfam{\vt_t}_{t\in\bS}$ be a strict \hl{\nbd{E_0}semi\-group} on $\sB^a(E)$, that is, a semigroup of unital endomorphisms $\vt_t$ of $\sB^a(E)$ that are strict. Meanwhile, there are several constructions of a product system from an \nbd{E_0}semigroup on $\sB^a(E)$; see \cite{Ske02,Ske03c,Ske05a,MSS06}. All these constructions capture, in a sense, the representation theory of $\sB^a(E)$. The first construction is due to Skeide \cite{Ske02}. This construction (inspired by Bhat's \cite{Bha96} for Hilbert spaces) is based on existence of a unit vector $\xi\in E$. The most general construction that works for arbitrary $E$ is based on the general representation theory of $\sB^a(E)$ in Muhly, Skeide and Solel \cite{MSS06}.

Let us discuss the construction based on \cite{MSS06}. We turn $E$ into a correspondence $_{\vt_t}E$ from $\sB^a(E)$ to $\cB$, by defining the left action $a.x:=\vt_t(a)x$. Since, by strictness of $\vt_t$, the action of the compacts on $_{\vt_t}E$ is nondegenerate, we may view $_{\vt_t}E$ as a correspondence from $\sK(E)$ to $\cB$. For every $t>0$ we define $E_t:=E^*\odot{_{\vt_t}E}$. Note that $E_t$ is a correspondence over $\cB$ that, likewise, may be viewed as correspondence over $\cB_E$. (The left action of $\cB_E$ is nondegenerate; see Section \ref{IspecialE}). Then
\beq{\label{repMe}
E\odot E_t
~=~
E\odot(E^*\odot{_{\vt_t}E})
~=~
(E\odot E^*)\odot{_{\vt_t}E}
~=~
\sK(E)\odot{_{\vt_t}E}
~=~
{_{\vt_t}E}
}\eeq
suggests that $E\odot E_t$ and $E$ are isomorphic as correspondences from $\sK(E)$ to $\cB$ but also as correspondences from $\sB^a(E)$ to $\cB$. That is, $a\odot\id_t$ should coincide with $\vt_t(a)$. In fact, interpreting all the identifications in the canonical way (see Sections \ref{IspecialE} and \ref{Itp}), we obtain an isomorphism $E\odot E_t\rightarrow E$ by setting
\beq{
x\odot(y^*\odot_tz)
~\longmapsto~
\vt(xy^*)z,
}\eeq
where we write $x^*\odot_ty$ in order to indicate that an elementary tensor $x^*\odot y$ is to be understood in $E^*\odot{_{\vt_t}}E$. We extend the definition to $t=0$ by putting $E_0=\cB$ and choosing the canonical identification $E\odot E_0=E$. (If $E$ is full, then this is automatic. Otherwise, we would find $E^*\odot{_{\vt_0}E}=E^*\odot E=\cB_E$.) The $E_t$ form a product system $E^\odot=\bfam{E_t}_{t\in\bS}$, that is
\beqn{
E_s\odot E_t
~=~
E_{s+t}
~~~~~~~~~~~~
(E_r\odot E_s)\odot E_t
~=~
E_r\odot(E_s\odot E_t),
}\eeqn
via
\bmun{
E_s\odot E_t
~=~
(E^*\odot{_{\vt_s}E})\odot(E^*\odot{_{\vt_t}E})
\\
~=~
E^*\odot{_{\vt_s}}(E\odot(E^*\odot{_{\vt_t}E}))
~=~
E^*\odot{_{\vt_s}}({_{\vt_t}}E)
~=~
E^*\odot{_{\vt_{s+t}}E}
~=~
E_{s+t}.
}\emun
We leave it as an instructive exercise to check on elementary tensors that the suggested identification
\beqn{
(x^*\odot_sy)\odot(x'^*\odot_ty')
~\longmapsto~
x^*\odot_{s+t}(\vt_t(yx'^*)y')
}\eeqn
is, indeed, associative.

We say the product system $E^\odot$ constructed before is the product system \hl{associated} with the \nbd{E_0}semigroup $\vt$. There are other ways to construct a product system of correspondences over $\cB$ from $\vt$, but they all lead to the same product system up to suitable isomorphism. (In the case of a von Neumann algebra $\cB$ there is the possibility to construct a product system of correspondences over the commutant $\cB'$; see Skeide \cite{Ske03c}. This product system is the \it{commutant} of all the others; see Section \ref{comndrsec}.) Our definition here is is for the sake of generality (it works for all strict \nbd{E_0}semigroups without conditions on $E$) and for the sake of uniqueness (it does not depend on certain choices like the choice of a unit vector in \cite{Ske02}).

Recall that for all $t>0$ the $E_t$ enjoy the property that they may also be viewed as correspondences over $\cB_E$. The uniqueness result \cite[Theorem 1.8 ]{MSS06} asserts that the $E_t$ are the only correspondences over $\cB_E$ that allow for an identification $E\odot E_t=E$ giving back $\vt_t(a)$ as $a\odot\id_t$. It is not difficult to show this statement remains true for the whole product system structure. We see also that the range ideal of $E_t$ cannot be smaller than $\cB_E$. Therefore, passing from $\cB$ to $\cB_E$ as \nbd{C^*}algebra, we may assume that $E^\odot$ is a \hl{full} product system, that is, that all $E_t$ $(t\in\bS)$ are full.

Now suppose we start with a full product system $E^\odot$. In order to establish that $E^\odot$ is (up isomorphism) the product system associated with a strict \nbd{E_0}semigroup, it is sufficient to find a full Hilbert module $E$ and identifications $E\odot E_t=E$ such that we have associativity
\beq{\label{asscond}
(E\odot E_s)\odot E_t
~=~
E\odot(E_s\odot E_t).
}\eeq
In that case, $\vt_t(a):=a\odot\id_t$ defines an \nbd{E_0}semigroup (Condition \eqref{asscond} gives the semigroup property) and the product system of this semigroup is
\beqn{
E^*\odot{_{\vt_t}E}
~=~
E^*\odot(E\odot E_t)
~=~
(E^*\odot E)\odot E_t
~=~
\cB\odot E_t
~=~
E_t.
}\eeqn
Suppose $E^\odot$ is a product system with a unital unit $\xi^\odot$. By a \hl{unit} for a product system $E^\odot$ we mean a family $\xi^\odot=\bfam{\xi_t}_{t\in\bS}$ of elements $\xi_t\in E_t$ with $\xi_0=\U$ that fulfills $\xi_s\odot\xi_t=\xi_{s+t}$. Note that this implies, in particular, that $\cB$ is unital. For nonunital $\cB$ we leave the term \it{unit} undefined! The unit is \hl{unital}, if all $\xi_t$ are unit vectors. (In particular, if $E^\odot$ has a unital unit, then $E^\odot$ is full.) It is well known that in this situation it is easy to construct an \nbd{E_0}semigroup. We merely sketch the construction and refer the reader to Bhat and Skeide \cite{BhSk00,Ske02} for details. For every $s,t\in\bS$ the map $\xi_s\odot\id_t\colon x_t\mapsto\xi_s\odot x_t$ defines an isometric embedding (as right module) of $E_t$ into $E_{s+t}$. The family of embeddings forms an inductive system, so that we may define the inductive limit $E_\infty=\lim_{t\to\infty}E_t$. For every $t\in\bS$ the factorization $E_s\odot E_t=E_{s+t}$ survives the inductive limit over $s$ and gives rise to a factorization $E_\infty\odot E_t=E_{\text{``$\infty+t$''}}=E_\infty$. Clearly, these factorizations fulfill \eqref{asscond}. Moreover, $E$ contains a unit vector, namely, the image $\xi$ of the vectors $\xi_t$ (which all coincide under the inductive limit). In particular, $E$ is full so that the product system of the \nbd{E_0}semigroup defined by setting $\vt_t(a):=a\odot\id_t$ is, indeed, $E^\odot$.

\bob \label{P1ob}
For Problem 1, which occupies the first half of these notes, this means the following: Suppose $E$ is a correspondence over $\cB$ with a unit vector $\xi$. Then $E^\odot=\bfam{E_n}_{n\in\N_0}$ with $E_n:=E^{\odot n}$ is a (discrete) product system and $\xi^\odot=\bfam{\xi_n}_{n\in\N_0}$ with $\xi_n:=\xi^{\odot n}$ is a unital unit. The inductive limit $E_\infty$ over that unit carries a strict \nbd{E_0}semigroup $\vt=\bfam{\vt_n}_{n\in\N_0}$ with $\vt_n(a)=a\odot\id_{E_n}$ whose product system is $E^\odot$. In particular, $E=E_1$ occurs as the correspondence of the unital strict endomorphism $\vt_1$ of $\sB^a(E_\infty)$.
\eob

We discuss briefly what the preceding construction does in the case of the \hl{trivial} product system ($E_n=\cB$ with product as left action and as product system operation) with a nontrivial unit vector (a proper isometry).

\bex\label{SzNFex}
Let $\cB$ denote a unital \nbd{C^*}algebra with a proper isometry $v\in\cB$. Then the inductive limit over the trival product system $\bfam{\cB^{\odot n}}_{n\in\N_0}$ with respect to the unit $\bfam{v^{\odot n}}_{n\in\N_0}$ has the form
\beq{\label{indds}
F
~:=~
\cB\oplus\bigoplus_{k=1}^\infty\cB_0
}\eeq
where $\cB_0:=(\U-vv^*)\cB$, and the induced endomorphism $\vt$ of $\sB^a(F)$ is $\vt(a)=uau^*$ where $u$ is the unitary defined by
\beqn{
u
~=~
v_0^*\oplus\id
\colon
\cB\oplus\bigoplus_{k=1}^\infty \cB_0
~\longrightarrow~
(\cB\oplus \cB_0)\oplus\bigoplus_{k=1}^\infty \cB_0
~=~
\cB\oplus\bigoplus_{k=0}^\infty \cB_0
~=~
\cB\oplus\bigoplus_{k=1}^\infty \cB_0,
}\eeqn
(in the last step we simply shift). It is an intriguing exercise to show that, indeed, the product system of $\vt$ is the trivial one (by \it{general abstract nonsense} this is true for every inner automorphism, but we mean to follow the construction from the beginning of this section), and to see how the embeddings $\cB=\cB^{\odot n}\rightarrow v^{\odot m}\odot\cB^{\odot n}=v^m\cB\subset\cB^{\odot(m+n)}=\cB$ really work and sit in $F$; see the old version \cite{Ske04p}.

In the case when $\cB=\sB(G)$ for some Hilbert space, we obtain just the Sz.-Nagy-Foias dilation of an isometry to a unitary.
\eex

In Sections \ref{C*uvsec} -- \ref{C*sec} it will be our job to reduce the cases we treat in these notes, full \nbd{C^*}modules over unital \nbd{C^*}algebras and strongly full \nbd{W^*}modules, to the case with a unit vector. We just mention that all results in the present section have analogues for \nbd{W^*}modules replacing strict mappings with normal (or $\sigma$--weak) mappings, replacing the tensor product of \nbd{C^*}cor\-re\-spond\-ences with that of \nbd{W^*}cor\-re\-spond\-ences, and replacing the word full by strongly full.

\section{Unit vectors in Hilbert modules}\label{C*uvsec}

In this section we discuss when full Hilbert modules over unital \nbd{C^*}algebras have unit vectors. In particular, we show that even if there is no unit vector, then a finite direct sum will admit a unit vector. This result will play its role in the solution of our Problem 1 in Theorem \ref{strunithm} for full correspondences over unital \nbd{C^*}algebras. As an application, not related to what follows, we give a simple proof of a statement about finitely generated Hilbert modules.

Of course, a Hilbert module $E$ over a unital \nbd{C^*}algebra $\cB$ that is not full cannot have unit vectors. But also if $E$ is full this does not necessarily imply existence of unit vectors.

\bex\label{non1ex}
Let $\cB=\C\oplus M_2=\tMatrix{\C&0\\0&M_2}\subset M_3=\sB(\C^3)$. The \nbd{M_2}\nbd{\C}module $\C^2=M_{21}$ may be viewed as a correspondence over $\cB$ (with operations inherited from $M_3\supset\tMatrix{0&0\\\C^2&0}$). Also its dual, the \nbd{\C}\nbd{M_2}module ${\C^2}^*=M_{12}=:\C_2$, may be viewed as a correspondence over $\cB$. It is easy to check that $M=\C^2\oplus\C_2=\tMatrix{0&\C_2\\\C^2&0}$ is a Morita equivalence (see Section \ref{MoritaSEC}) from $\cB$ to $\cB$ (in particular, $M$ is full) without a unit vector.

Note that $M\odot M=\cB$ has a unit vector. Example \ref{lastEX} tells us that there are serious examples in the discrete case where not one of the tensor powers $E^{\odot n}$ $(n>0)$ has a unit vector.
\eex

Observe that all modules and correspondences in Example \ref{non1ex} are \nbd{W^*}modules, so missing unit vectors are not caused by insufficient closure. The reason why $M$ does not contain a unit vector is because the full Hilbert \nbd{M_2}module $\C_2$ has ``not enough space'' to allow for sufficiently many orthogonal vectors. (Not two nonzero vectors of this module are orthogonal.) Another way to argue is to observe that every nonzero inner product $\AB{x^*,y^*}$ is a rank-one operator in $M_2=\sB(\C^2)$ while the identity has rank two. As soon as we create ``enough space'', for instance, by taking the direct sum of sufficiently many (in our case two) copies of $\C_2$ the problem disappears.

In the following lemma we show that for every full Hilbert module a finite number of copies will be ``enough space''. The basic idea is that, if $\AB{x,y}=\U$, then by Cauchy-Schwartz inequality $\U=\AB{x,y}\AB{y,x}\le\AB{x,x}\norm{y}^2$ so that $\AB{x,x}$ is invertible and $x\sqrt{\AB{x,x}^{-1}}$ is a unit vector. Technically, the condition $\AB{x,y}=\U$ is realized only approximately and by elements in $E^n$ rather than in $E$.

\blem\label{univeclem}
Let $E$ be a full Hilbert module over a unital \nbd{C^*}algebra. Then there exists $n\in\N$ such that $E^n$ has a unit vector.
\elem

\proof
$E$ is full, so there exist $x_i^n,y_i^n\in E$ $(n\in\N;i=1,\ldots,n)$ such that
\beqn{
\lim_{n\to\infty}\sum_{i=1}^n\AB{x_i^n,y_i^n}
~=~
\U.
}\eeqn
The subset of invertible elements in $\cB$ is open. Therefore, for $n$ sufficiently big $\sum_{i=1}^n\AB{x_i^n,y_i^n}$ is invertible. Defining the elements $X_n=(x_1^n,\ldots,x_n^n)$ and $Y_n=(y_1^n,\ldots,y_n^n)$ in $E^n$ we have, thus, that
\beqn{
\AB{X_n,Y_n}
~=~
\sum_{i=1}^n\AB{x_i^n,y_i^n}
}\eeqn
is invertible. So, also $\AB{X_n,Y_n}\AB{Y_n,X_n}$ is invertible and, therefore, bounded below by a strictly positive constant. Of course, $\norm{Y_n}\ne0$. By Cauchy-Schwartz inequality also
\beqn{
\AB{X_n,X_n}
~\ge~
\frac{\AB{X_n,Y_n}\AB{Y_n,X_n}}{\norm{Y_n}^2}
}\eeqn
is bounded below by a strictly positive constant and, therefore, $\AB{X_n,X_n}$ is invertible. It follows that $X_n\sqrt{\AB{X_n,X_n}^{-1}}$ is a unit vector in $E^n$.\qed

\bcor\label{univeccor}
If $E$ (as before) contains an arbitrary number of mutually orthogonal copies of a full Hilbert submodule (for instance, if $E$ is isomorphic to $E^n$ for some $n\ge2$), then $E$ has a unit vector.
\ecor

Lemma \ref{univeclem} implies that, if $\sK(E)$ is unital, then $\sK(E)=\sF(E)$. (Just apply the lemma to the full Hilbert \nbd{\sK(E)}module $E^*$.)

\bcor\label{Kunicor}
If $\sK(E)$ is unital, then $E$ is algebraically finitely generated.
\ecor

This is some sort of inverse to the well-known fact that an (algebraically) finitely generated Hilbert \nbd{\cB}module is isomorphic to a (complemented) submodule of $\cB^n$ for some $n$.

\section{Unit vectors in $W^*$--modules}\label{W*uvsec}

In this section we proof the analogue of Lemma \ref{univeclem} for \nbd{W^*}modules. Of course, a \nbd{W^*}module is a Hilbert module. If it is full then Lemma \ref{univeclem} applies. But the good notion of fullness for a \nbd{W^*}module is that it is \hl{strongly full}, that is, the inner product of the \nbd{W^*}module generates $\cB$ as a \nbd{W^*}algebra. (Strong fullness is the more useful notion for \nbd{W^*}modules, because it can always be achieved by restricting $\cB$ to the \nbd{W^*}subalgebra generated by the inner product. Example \ref{sfnfex} tells us that the same is not true for fullness in the case of \nbd{W^*}modules.) It is the assumption of strong fullness for which we want to resolve Problem 1 for \nbd{W^*}modules, and not the stronger assumption of fullness (that might be not achievable). We thank B.\ Solel for pointing out to us this gap in the first version of these notes.

We see immediately that for strongly full \nbd{W^*}modules the cardinality of the direct sum in Lemma \ref{univeclem} can no longer be kept finite.

\bex\label{sfnfex}
Let $H$ be an infinite-dimensional Hilbert space. Then $H^*$ is a \nbd{W^*}module over $\sB(H)$, that is strongly full but not full as a Hilbert \nbd{\sB(H)}module. (Indeed, the range ideal of $H^*$ in $\sB(H)$ is $\sB(H)_{H^*}=\sK(H)\ne\sB(H)$.) For every finite direct sum ${H^*}^n$ the inner product $\AB{X_n,X_n}$ $(X_n\in{H^*}^n)$ has rank not higher than $n$. Therefore, ${H^*}^n$ does not admit a unit vector. Only if we consider $\ol{{H^*}^\en}^s$, the self-dual extension of ${H^*}^\en$, where $\en=\dim H$, then the vector in $\ol{{H^*}^\en}^s$ with the components $e_i^*$ ($\bfam{e_i}$ some orthonormal basis of $H$) is a unit vector. But this vector is not in ${H^*}^\en$ if $\en$ is infinite.

Observe that, for arbitrary cardinality $\en$, we have ${H^*}^\en=\sK(H,\C^\en)$, while $\ol{{H^*}^\en}^s=\sB(H,\C^\en)$. In fact, when $\dim H=\en$ we have $\ol{{H^*}^\en}^s=\sB(H)$.
\eex

The example is in some sense typical. In fact, we constructed a multiple of $H^*$ that contains a unit vector by choosing an orthonormal basis for its dual $H$. This will also be our strategy for general \nbd{W^*}modules. A suitable substitute for orthonormal bases are \it{quasi} orthonormal bases. A \hl{quasi orthonormal basis} in a \nbd{W^*}module $E$ over $\cB$ is a family $\bfam{e_i,p_i}_{i\in S}$ where $S$ is some index set (of cardinality $\en$, say), $p_i$ are projections in $\cB$ and $e_i$ are elements in $E$ such that
\beqn{
\AB{e_i,e_j}
~=~
\delta_{i,j}p_j
\text{~~~~~~and~~~~~~}
\sum_{i\in S}e_ie_i^*
~=~
\id_E
}\eeqn
(monotone limit in the \nbd{W^*}algebra $\sB^a(E)$ over the finite subsets of $S$ in the case $S$ is not finite). Existence of a quasi orthonormal basis follows from self-duality of $E$ and monotone completeness of $\sB^a(E)$ by an application of \it{Zorn's lemma}; see Paschke \cite{Pas73}.

\blem\label{Wuniveclem}
Let $E$ be a strongly full \nbd{W^*}module. Then there exists a cardinal number $\en$ such that $\ol{E^\en}^s$ has a unit vector.
\elem

\proof
Let us choose a quasi orthonormal basis $\bfam{e_i^*,e_ie_i^*}_{i\in S}$ for the dual \nbd{\sB^a(E)}module $E^*$. (Observe that $E^*$ is a \nbd{W^*}module; see Remark \ref{dctrem}.) Then
\beqn{
\sum_{i\in S}e_i^*e_i
~=~
\sum_{i\in S}\AB{e_i,e_i}
~=~
\id_E.
}\eeqn
The second sum is, actually, over the elements $\AB{e_i,e_i}$ when considered as operator acting from the left on $E^*$. But, as $E$ is strongly full, the action of $\cB$ on $E^*$ is faithful. In particular, the only element in $\cB$ having the action $\id_E$ is, really, $\U\in\cB$. Now, if we put $\en=\#S$, then the vector in $\ol{E^\en}^s$ with components $e_i$ is a unit vector.\qed

\lf
In Gohm and Skeide \cite{GoSk05} we pointed out that existence of a quasi orthonormal basis for a \nbd{W^*}module may be used to give a simple proof of the \it{amplification-induction theorem}, that is, the theory of normal representations of a von Neumann algebra $\cB$. Indeed, let $\cB\subset\sB(G)$ be a von Neumann algebra acting nondegenerately on a Hilbert space $G$. If $\rho$ is a nondegenerate representation of $\cB$ on another Hilbert space $H$. Then $E':=\bCB{x'\in\sB(G,H)\colon \rho(b)x'=x'b~(b\in\cB)}$ is a \nbd{W^*}module over $\cB'\subset\sB(G)$ with inner product $\AB{x',y'}:=y'^*x'\in\cB'$. Moreover, $\cls E'G=H$; see Section \ref{comndrsec}, in particular, Remark \ref{dctrem}. Let $\bfam{e'_i,p'_i}_{i\in S}$ be a quasi orthonormal basis of $E'$. It follows that $H=\bigoplus_{i\in S}p'_iG\subset G^{\#S}=G\otimes\C^{\#S}$ (see Section \ref{Icard} for notation). The representation $\rho$ is, then, the compression of the amplification $\id_\cB\otimes\id_{\C^{\#S}}$ to the invariant subspace $H$.

We may use Lemma \ref{Wuniveclem} to furnish a new proof of the structure theorem for algebraic isomorphisms of von Neumann algebras. Indeed, let $\rho$ be faithful so that (see Section \ref{comndrsec}) $E'$ is strongly full. By Lemma \ref{Wuniveclem} a suitable multiple $\ol{E'^\en}^s$ of $E'$ contains a unit vector $\xi'$. We may choose a quasi orthonormal basis $\bCB{(\xi',\U)}\cup\bfam{e'_i,p'_i}_{i\in S}$ of $\ol{E'^\en}^s$ (disjoint union). Let $\el$ be the smallest infinite cardinal number not smaller than $\#S$. Then the multiple $\ol{E'^{\en\cdot\el}}^s$ of $\ol{E'^\en}^s$ is isomorphic to $\soplus_{i\in S}(E'_i)^\el$ where $E'_i:=\cB'\oplus p'_i\cB'=(\U-p'_i)\cB'\oplus p'_i\cB'\oplus p'_i\cB'$. It follows that $\ol{E'^\el_i}^s\cong\ol{(\U-p'_i)\cB'^\el}^s\oplus\ol{p'_i\cB'^\el}^s=\ol{\cB'^\el}^s$. In other words, $\ol{E'^{\en\cdot\el}}^s\cong\ol{\cB'^\el}^s$.

\bcor\label{ampisocor}
If $\rho$ is a faithful normal nondegenerate representation of a von Neumann algebra $\cB\subset\sB(G)$ on $H$, then there exists a Hilbert space $\eH$ such that the representations $b\mapsto\rho(b)\otimes\id_\eH$ and $b\mapsto b\otimes\id_\eH$ are unitarily equivalent.
\ecor

\section{Morita equivalence for product systems}\label{MoritaSEC}

In this section we review the notions of (strong) Morita equivalence (Rieffel \cite{Rie74a}), Morita equivalence for Hilbert modules (new in these notes) and Morita equivalence for correspondences (Muhly and Solel \cite{MuSo00}). We put some emphasis on the difference between the \nbd{C^*}case and the \nbd{W^*}case. That difference is in part responsible for the fact that we can solve Problem 1 in full generality only for \nbd{W^*}modules. The \nbd{C^*}case can be done only for unital \nbd{C^*}algebras and, even under this assumption, it is much less elegant. Then we show that a product system of \nbd{W^*}correspondences can be derived from an \nbd{E_0}semigroup, if and only if it is Morita equivalent to a product system that has a unital unit. In the discrete case this means a \nbd{W^*}correspondence stems form a unital endomorphism of some $\sB^a(E)$, if and only if it is Morita equivalent to a \nbd{W^*}correspondence that has a unit vector.

A correspondence $M$ from $\cA$ to $\cB$ is called a \hl{Morita equivalence} from $\cA$ to $\cB$, if it is full and if the canonical mapping from $\cA$ into $\sB^a(M)$ corestricts to an isomorphism $\cA\rightarrow\sK(M)$. Clearly, the two conditions can be written also as
\baln{
M^*\odot M
&
~=~
\cB
&
M\odot M^*
&
~=~
\cA.
}\ealn
From these equations one concludes easily a couple of facts. Firstly, if $M$ is a Morita equivalence from $\cA$ to $\cB$, then $M^*$ is a Morita equivalence from $\cB$ to $\cA$. Secondly, the tensor product of Morita equivalences is a Morita equivalence. Thirdly, $M$ and $M^*$ are inverses under tensor product. Two \nbd{C^*}algebras are called \hl{strongly Morita equivalent}, if they admit a Morita equivalence from one to the other. Usually, we say just \it{Morita equivalent} also when we intend \it{strongly Morita equivalent}.

\bex\label{Mnex}
All $M_n$ are Morita equivalent to $\C$ via the Morita equivalence $\C^n$.
\eex

\bex\label{repthex}
Also the representation theory of $\sB^a(E)$ is just a matter of Morita equivalence. In fact, the identity of the \nbd{\sK(E)}\nbd{\cB}correspondences in \eqref{repMe} becomes \it{crystal}, taking into account that $E$ is a Morita equivalence from $\sK(E)$ to $\cB_E$ and $E^*$ is its inverse; see \cite{MSS06}.
\eex

In the category of \nbd{W^*}algebras with \nbd{W^*}correspondences a correspondence from $\cA$ to $\cB$ is a \hl{Morita \nbd{W^*}equivalence}, if $M$ is strongly full and if the canonical mapping $\cA\rightarrow\sB^a(M)$ is an isomorphism.

\brem
Clearly, in the \nbd{W^*}case we have $M\sodots M^*=\sB^a(M)$. The fact that Morita equivalence for \nbd{W^*}algebra relates $\cA$ to $\sB^a(M)$ while strong Morita equivalence of \nbd{C^*}algebras relates $\cA$ only to $\sK(M)$ is one of the reasons why our solution of Problem 1 works only in the \nbd{W^*}case, respectively, runs considerably less smoothly in the particular \nbd{C^*}case we discuss in Section \ref{C*sec}.
\erem

\bex\label{WMnex}
The $M_n$ are \nbd{W^*}algebras, the $\C^n$ and their duals are \nbd{W^*}correspondences and all tensor products are tensor products in the \nbd{W^*}sense. So, Example \ref{Mnex} is also an example for Morita equivalence of \nbd{W^*}algebras.
\eex

\brem
Versions of Examples \ref{Mnex} and \ref{WMnex} for infinite-dimensional matrices and $\C$ replaced with $\cB$ are crucial to resolve Problem 1. Essentially, we are going to use $\cB^\en$ as Morita equivalence from $M_\en(\cB)$ to $\cB$. Of course, for infinite-dimensional matrices either we have to pass to strong closures (Section \ref{mthmsec}) or to a weaker notion of Morita equivalence (Section \ref{C*sec}).
\erem

\bemp[Definition (Muhly and Solel \cite{MuSo00}).~]
A correspondence $E$ over $\cB$ and a correspondence $F$ over $\cC$ are \hl{Morita equivalent}, if there is a Morita equivalence $M$ from $\cB$ to $\cC$ such that $E\odot M=M\odot F$ (or $E=M\odot F\odot M^*$).
\eemp

We add here:

\bdefi
A Hilbert \nbd{\cB}module $E$ and a Hilbert \nbd{\cC}module $F$ are \hl{Morita equivalent}, if there is a Morita equivalence $M$ from $\cB$ to $\cC$ such that $E\odot M=F$ (or $E=F\odot M^*$).
\edefi

Of course, the definitions for the \nbd{W^*}case are analogue.

Morita equivalence of Hilbert modules and Morita equivalence of correspondences are related by the following crucial proposition. Suppose $\alpha\colon\sB^a(E)\rightarrow\sB^a(F)$ is a (bi-)strict isomorphism. By \cite{MSS06} this is the case, if and only if $E$ and $F$ are Morita equivalent where the Morita equivalence $M$ induces $\alpha$ as $\alpha(a)=a\odot\id_M$.

Now suppose there are two strict unital endomorphisms $\vt$ and $\theta$ on $\sB^a(E)$ and $\sB^a(F)$, respectively. We may ask whether they are \hl{conjugate}, that is, whether there exists a (bi-)strict isomorphism $\alpha\colon\sB^a(E)\rightarrow\sB^a(F)$ such that $\theta=\alpha\circ\vt\circ\alpha^{-1}$.

\bprop\label{conprop}
$\vt$ and $\theta$ are conjugate, if and only if there is a Morita equivalence inducing an isomorphism $F=E\odot M$ such that $E_\vt\odot M=M\odot F_\theta$, that is, if and only if $E$ and $F$ as well as $E_\vt$ and $F_\theta$ are Morita equivalent by the same Morita equivalence. 
\eprop

The proof consists very much of computations like the second half of the proof Theorem \ref{psMethm} below. We leave it as an exercise.

\brem
Note that in the scalar case $\cB=\cC=\C$, where $\C$ is the only Morita equivalence over $\C$, we recover the well-known facts that every normal isomorphism $\alpha\colon\sB(G)\rightarrow\sB(H)$ is induced by a unitary $G\rightarrow H$ and that the multiplicity spaces of two endomorphisms conjugate by $\alpha$ must be equal.
\erem

Clearly, if $E^\odot=\bfam{E_t}_{t\in\bS}$ is a product system of correspondences over $\cB$ and $M$ is a Morita equivalence from $\cB$ to $\cC$, then $M^*\odot E^\odot\odot M:=F^\odot=\bfam{F_t}_{t\in\bS}$  with $F_t:=M^*\odot E_t\odot M$ and isomorphisms
\beqn{
F_s\odot F_t
~=~
M^*\odot E_s\odot M\odot M^*\odot E_t\odot M
~=~
M^*\odot E_s\odot E_t\odot M
~=~
M^*\odot E_{s+t}\odot M
~=~
F_{s+t}
}\eeqn
is a product system of correspondences over $\cC$.

\bdefi
We say $E^\odot$ and $F^\odot$ are \hl{Morita equivalent}, if there exists a Morita equivalence $M$ and an \hl{isomorphism} $u^\odot=\bfam{u_t}_{t\in\bS}\colon M^*\odot E^\odot\odot M\rightarrow F^\odot$ (that is, the $u_t$ are bilinear unitaries $M^*\odot E_t\odot M\rightarrow F_t$ such that $u_s\odot u_t=u_{s+t}$ and $u_0=\id_\cC$).
\edefi

The version for \nbd{W^*}correspondences is analogue.

The following corollary is proved very much like Proposition \ref{conprop} taking also into account (see \cite{Ske02}) that two strict \nbd{E_0}semigroups on the same $\sB^a(E)$ are cocycle conjugate, if and only if their product systems are isomorphic.

\bcor\label{concor}
Suppose $\vt$ and $\theta$ are strict \nbd{E_0}semigroups on $\sB^a(E)$ and on $\sB^a(F)$, respectively. Then $\vt$ and $\theta$ are cocycle conjugate via a (bi-)strict isomorphism $\alpha\colon\sB^a(E)\rightarrow\sB^a(F)$ (in the sense that $\theta$ and $\alpha\circ\vt\circ\alpha^{-1}:=\bfam{\alpha\circ\vt_t\circ\alpha^{-1}}_{t\in\bS}$ are cocycle equivalent), if and only if the product systems $E^\odot$ of $\vt$ and $F^\odot$ of $\theta$ are Morita equivalent via the Morita equivalence $M$ that induces $\alpha$ as $\alpha(a)=a\odot\id_M$.
\ecor

Of course, also here there is a version for \nbd{W^*}modules.

\bthm\label{psMethm}
Let $E^{\sodots}=\bfam{E_t}_{t\in\bS}$ be a product system of strongly full \nbd{W^*}correspondences $E_t$ over a \nbd{W^*}al\-ge\-bra $\cB$. Then $E^{\sodots}$ is the product system of a normal \nbd{E_0}semigroup $\vt=\bfam{\vt_t}_{t\in\bS}$ on $\sB^a(E)$ for some \nbd{W^*}module $E$ over $\cB$, if and only if $E^{\sodots}$ is Morita equivalent to a product system $F^{\sodots}$ of \nbd{W^*}correspondences over a \nbd{W^*}algebra $\cC$ that contains a unital unit $\zeta^\odot$.
\ethm

\proof
``$\Longrightarrow$''. Suppose $E^{\sodots}$ is the product system of the normal \nbd{E_0}semigroup $\vt$ on the \nbd{W^*}algebra $\cC:=\sB^a(E)$. Put $F_t:=E\sodots E_t\sodots E^*$. As $E\sodots E_t=E$ and $\vt_t(a)=a\odot\id_{E_t}$, we find $F_t={_{\vt_t}}\sB^a(E)$ and $a_s\odot a_t=\vt_t(a_s)a_t$ is the isomorphism $F_s\sodots F_t=F_{s+t}$. Clearly, $\zeta_t=\id_{E}\in\sB^a(E)=F_t$ defines a unital unit $\zeta^\odot$ for $F^{\sodots}$. (One easily verifies that also the inductive limit $F_\infty=\sB^a(E)=E\sodots E^*$ constructed from that unit is that obtained from $E$ via the Morita equivalence $M:=E^*$ as $F_\infty=E\sodots E^*$.)

``$\Longleftarrow$''. Suppose $M$ is a Morita \nbd{W^*}equivalence from $\cB$ to $\cC$ such that $F^{\sodots}:=M^*\sodots E^{\sodots}\sodots M$ has a unital unit $\zeta^\odot$. Construct the inductive limit $F_\infty$ with the normal \nbd{E_0}semigroup $\theta_t(a)=a\odot\id_{F_t}$ $(a\in\sB^a(F_\infty))$ and put $E:=F_\infty\sodots M^*$. Then
\beqn{
\vt_t(a)
~:=~
\theta_t(a\odot\id_M)\odot\id_{M^*}
~=~
a\odot\id_M\odot\id_{F_t}\odot\id_{M^*}
~=~
a\odot\id_{E_t}
}\eeqn
$(a\in\sB^a(E))$ where
\baln{
E
~=~
E\sodots M\sodots M^*
&~(=~
F_\infty\sodots M^*
~=~
F_\infty\sodots F_t\sodots M^*)
\\
&~=~
E\sodots M\sodots F_t\sodots M^*=E\sodots E_t
}\ealn
is the induced semigroup on $\sB^a(E)$. As $F_\infty$ is full (it contains a unit vector) also $E=F_\infty\sodots M^*$ is full ($\ol{\cC_{F_\infty}}^s=\cC$ acts nondegenerately on $M^*$ so that $\ol{\cB_E}^s=\ol{\cB_{M^*}}^s=\cB$). Therefore, by uniqueness the product systems associated with $\vt$ gives us back $E_t$.\qed

\bcor
Two product systems of \nbd{W^*}correspondences in the same Morita equivalence class are either both or are both not the product systems of a normal \nbd{E_0}semigroups.
\ecor

\section{Endomorphisms: \nbd{W^*}case}\label{mthmsec}

The results of Sections \ref{W*uvsec} and \ref{MoritaSEC} allow in a very plain way to resolve Problem 1 for \nbd{W^*}cor\-re\-spond\-ences. But, we explain first the idea in the \nbd{C^*}case under the assumptions of Lemma \ref{univeclem} --- and in doing so we illustrate why it does not work in the \nbd{C^*}case. This helps to appreciate better the \nbd{W^*}case.

Let $E$ be a full correspondence over a unital \nbd{C^*}algebra $\cB$. By Lemma \ref{univeclem} we know that for some $n\in\N$ the correspondence $E^n$ has a unit vector. We observe that $E^n=\cB^n\odot E$, where $\cB^n$ is a Morita equivalence from $M_n(\cB)$ to $\cB$. If we could show existence of a unit vector in $M_n(E)=\cB^n\odot E\odot\cB_n$ where $\cB_n:=(\cB^n)^*$ is the dual of $\cB^n$, then $E$ was Morita equivalent to a correspondence with a unit vector. In this case the ``$\Longleftarrow$'' direction of the proof of Theorem \ref{psMethm} works even without strong closure. (One main reason for strong closure is that rarely $\sB^a(E)=\sK(E)$ so $E$ is a rarely a Morita equivalence from $\sB^a(E)$ to $\cB_E$ as needed in the proof of Theorem \ref{psMethm}. But, here with $\cB$ also $M_n(\cB)=\sB^a(\cB^n)=\sK(\cB^n)$ is unital.)

Unfortunately, $M_n(E)$ need not have a unit vector. Suppose $n\ge2$ is the minimal cardinality such that $E^n$ has a unit vector. To produce a $\U$ in a place in the diagonal we need $n$ orthogonal vectors, and to produce $\U$ in each of the $n$ places in the diagonal we need $n^2$ orthogonal vectors. However, the \nbd{M_n(\cB)}correspondence $M_n(E)$ still has ``space'' only for $n$ orthogonal vectors with suitable inner products. We invite the reader to check that for the correspondence $E$ from Example \ref{non1ex}, where $n=2$, the correspondence $M_2(E)$ does admit unit vectors. The problem remains, when we use $M_m(E)$ $(m>n)$ instead. It disappears if $m=\infty$, because then we can ``slice'' $m=\infty$ into $n$ slices still of size $m=\infty$. The problem is now that the sums when calculating inner products of elements in $M_\infty(E)$ (or also in products of elements in $M_\infty(\cB)$) converge only strongly. (For instance, $\U_\infty\in M_\infty(\cB)$ is approximated by $\U_m\in M_m(\cB)\subset M_\infty(\cB)$.) This is a second reason why we have to switch to the \nbd{W^*}case.

In the context of \nbd{W^*}modules, Lemma \ref{Wuniveclem} allows for arbitrary cardinalities $\en$. We start by giving a precise meaning to $M_\en(\cB)$ and $M_\en(E)$. So let $E$ be a \nbd{W^*}cor\-re\-spond\-ences over a \nbd{W^*}algebra $\cB$. Let $S$ be a set with cardinality $\#S=\en$ and denote by $\bfam{e_k}_{k\in S}$ the natural orthonormal basis of $\C^\en$. We set $M_\en(\cB):=\sB(\C^\en)\sbars{\otimes}\cB$ (tensor product of \nbd{W^*}algebras) and we identify an element $B\in M_\en(\cB)$ with the matrix $\bfam{b_{ij}}_{i,j\in S}$ where
\beqn{
b_{ij}
~=~
(e_i^*\otimes\id_\cB)B(e_j\otimes\id_\cB)
~\in~
\cB.
}\eeqn
We put $M_\en(E):=\sB(\C^\en)\sbars{\otimes}E$, that is, the \hl{exterior tensor product} of \nbd{W^*}modules; see \cite[Section 4.3]{Ske01}. We identify an element $X\in M_\en(E)$ with the matrix $\bfam{x_{ij}}_{i,j\in S}$ where
\beqn{
x_{ij}
~=~
(e_i^*\otimes\id_E)X(e_j\otimes\id_E)
~\in~
E.
}\eeqn
The operations in this correspondence over $M_\en(\cB)$ are
\bal{\label{matop}
\AB{X,Y}_{ij}
~&=~
\sum_k\AB{x_{ki},y_{kj}}
&
(XB)_{ij}
~&=~
\sum_kx_{ik}b_{kj}
&
(BX)_{ij}
~&=~
\sum_kb_{ik}x_{kj},
}\eal
where all sums are \nbd{\sigma}strong limits. A matrix $X=\bfam{x_{ij}}$ is an element of $M_\en(E)$, if and only if all $\sum_k\AB{x_{ki},x_{kj}}$ exist \nbd{\sigma}strongly and define the matrix elements of an element in $M_\en(\cB)$.

Clearly, $M:=\C^\en\sbars{\otimes}\cB=\ol{\cB^\en}^s$ is a Morita \nbd{W^*}equivalence from $M_\en(\cB)$ to $\cB$ and
\beqn{
M\sodots E\sodots M^*
~=~
M_\en(E).
}\eeqn

\bcor\label{inftycor}
$E$ and $M_\en(E)$ are Morita equivalent \nbd{W^*}correspondences.
\ecor

\bprop\label{inftyprop}
Suppose $E$ is strongly full and let $\en$ be an \hl{infinite} cardinal number not smaller than that granted by Lemma \ref{Wuniveclem}. Then $M_\en(E)$ has a unit vector.
\eprop

\proof
Denote by $\el$ the cardinal number from Lemma \ref{Wuniveclem} and fix $\en$ as stated. Choose sets $S,T$ with $\#S=\el,\#T=\en$. Let $x_\ell$ $(\ell\in S)$ denote the components of a unit vector in $\ol{E^\el}^s$. As $\en$ is infinite (by assumption!) and $\el\le\en$ so that $\el\en=\en$, we may fix a bijection $\vp\colon T\rightarrow S\times T$. Denote by $\vp_1$ and $\vp_2$ the first and the second component, respectively, of $\vp$. Define a matrix $X\in M_\en(E)$ by setting
\beqn{
x_{ij}
~=~
x_{\vp_1(i)}\delta_{\vp_2(i),j}.
}\eeqn
Then
\bmun{
\AB{X,X}_{ij}
~=~
\sum_{k\in T}\AB{x_{ki},x_{kj}}
~=~
\sum_{k\in T}\delta_{\vp_2(k),i}\delta_{\vp_2(k),j}\AB{x_{\vp_1(k)},x_{\vp_1(k)}}
\\
~=~
\sum_{(\ell,k)\in S\times T}\delta_{k,i}\delta_{k,j}\AB{x_\ell,x_\ell}
~=~
\Bfam{\sum_{k\in T}\delta_{k,i}\delta_{k,j}}\Bfam{\sum_{\ell\in S}\AB{x_\ell,x_\ell}}
~=~
\delta_{i,j}\U.~~~\qedsymbol
}\emun
\noqed

\noindent
Putting together Corollary \ref{inftycor} and Proposition \ref{inftyprop} with Theorem \ref{psMethm} resolves Problem 1 in the \nbd{W^*}case.

\bthm\label{normunithm}
Let $E$ be a strongly full \nbd{W^*}correspondence. Then there is a \nbd{W^*}module $F$ (necessarily strongly full) and a unital normal endomorphism of $\vt$ of $\sB^a(F)$ such that $F_\vt=E$.
\ethm

\bob\label{uniexob}
Note that, by construction of $M_\en(E)$ and the proof of Theorem \ref{psMethm}, the module $F=\bfam{\limind_{n\to\infty}M_\en(E^{\sodots n})}\sodots\ol{\cB^\en}^s\supset M_\en(E)\sodots\ol{\cB^\en}^s=\ol{E^\en}^s$ contains a unit vector. But we can show even more. In Corollary \ref{freecor} we will see that, for a suitable cardinality $\en$, we may even achieve that $F$ is isomorphic to a \it{free} module $\ol{\cB^\en}^s$. This observation (and its \nbd{C^*}version, Observation \ref{C*uniexob}, below) are in duality with Observation \ref{repampob} in the sense of \it{commutant} (Section \ref{comndrsec}).
\eob

\brem\label{nsfrem}
If $E$ is not necessarily strongly full, then, as explained in Section \ref{psunisec} (\nbd{W^*}version, of course), for that $E=F_\vt$ for some $\vt$ acting on some $\sB^a(F)$, it is necessary that $\ol{\cB_E}^s$ acts nondegenerately on $E$. But this is also sufficient, for in this case, we apply Theorem \ref{normunithm} to the strongly full correspondence $E$ over $\ol{\cB_E}^s$ and obtain a $\vt$ acting on $\sB^a(F)$ for some $F$ that is strongly full over $\ol{\cB_E}^s$.
\erem

\section{Endomorphisms: $C^*$--case}\label{C*sec}

In this section we resolve Problem 1 for full \nbd{C^*}correspondences over a unital \nbd{C^*}al\-ge\-bra. The proof is less streamlined than that of the \nbd{W^*}case, so we do not develop a complete analogue of the treatment of the \nbd{W^*}version --- also because, partly, this is not possible.

One problem was to have a notion of Morita equivalence that understands a full Hilbert \nbd{\cB}module $E$ as a Morita equivalence from $\sB^a(E)$ to $\cB$ and not just from $\sK(E)$ to $\cB$. In the previous sections the strongly closed versions for \nbd{W^*}objects did the job. In this section we elaborate a version for \hl{strict} closure (or what is the same for strict or \nbd{*}strong completion). And we elaborate this \it{strict} Morita equivalence only for the case, where one of the algebras is $\sB^a(\cB)$. This will allow for the necessary matrix constructions, and Lemma \ref{univeclem} will guarantee existence of a unit vector in the matrix modules. The fact that, for nonunital \nbd{C^*}algebras, we have available neither Lemma \ref{univeclem} nor Lemma \ref{Wuniveclem} is responsible for that we cannot prove the result in that case. Lemma \ref{univeclem} works only for full Hilbert modules over unital \nbd{C^*}algebras, and the strict completion will be only ``strictly full'' over the multiplier algebra of $\cB$. The proof of Lemma \ref{Wuniveclem} is based on quasi orthonormal bases that, in strict completions, are not available.

\bex
Let $\cB=C_0(-2,2)$ and $I=C_0(-1,1)\subset\cB$ an ideal and define the full Hilbert \nbd{\cB}module $E=\cB\oplus I$. Then the strict completion of $E$, $\sB^a(\cB,E)$, is the direct sum $C_b(-2,2)\oplus C_0(-1,1)$. The only nonzero projection in $\sB^a(\cB)=C_b(-2,2)$ is $\U$. Every element $\xi$ in $\sB^a(\cB,E)$ that has unit length leaves a nonzero complement $(\U-\xi\xi^*)E$ and all inner products of elements in that complement are in $I$. So $\sB^a(\cB,E)$ has no quasi orthonormal basis.

Nevertheless, we remark that $E$, equipped with its natural left action, is the correspondence of a strict unital endomorphism on $\sB^a(F)$ for some full Hilbert \nbd{\cB}module $F$. Indeed, choose two Hilbert spaces $H_1$ and $H_2$ and put $F=C_0((-1,1),H_1)\oplus C_0((-2,2),H_2)$. Then $F\odot E=C_0((-1,1),H_1\oplus H_1\oplus H_2)\oplus C_0((-2,2),H_2)$. Therefore, whenever $H_1$ is infinite-dimensional, there exists an isomorphism $H_1\oplus H_1\oplus H_2\rightarrow H_1$ so that $F\odot E$ and $F$ are isomorphic. We do not know a counter example for nonunital $\cB$.
\eex

Let us start with some generalities, however, without discussing (as would be natural) how the definitions fit into the frame of multiplier algebras, double centralizers and strict topology. If $E$ is a Hilbert \nbd{\cB}module, then by the \hl{strict completion} of $E$ we understand the space $\sB^a(\cB,E)$. If $\cB$ is unital, then $\sB^a(\cB,E)$ is just $E$ where we consider $x\in E$ as the map $b\mapsto xb$. In general, $\sB^a(\cB,E)$ is Hilbert module over $\sB^a(\cB)$ with inner product $\AB{X,X'}=X^*X'$. Further, $\sB^a(\cB,E)$ has the same operators as $E$, that is, $\sB^a(\sB^a(\cB,E))=\sB^a(E)$. (An element $a\in\sB^a(E)$ acts on $X\in\sB^a(\cB,E)$ simply by composition $aX$, while an operator $a$ on $\sB^a(\cB,E)$ determines the operator $Xb\mapsto (aX)b$ on $E$.)

Now we wish to define an appropriate tensor product among such spaces.

\bprop\label{strtpprop}
Let $E$ be a Hilbert \nbd{\cB}module and let $F$ be a correspondence from $\cB$ to $\cC$. Then:
\begin{enumerate}
\item\label{str1}
The left action of $\cB$ on $F$ extends to a (unique and strict) action of $\sB^a(\cB)$. Therefore, also $\sB^a(\cC,F)$ has a left action of $\sB^a(\cB)$.

\item\label{str2}
For every $X\in\sB^a(\cB,E)$ by setting $\eta(X)=X\odot\id_F$ we define a map in $\sB^a(\cB\odot F,E\odot F)=\sB^a(F,E\odot F)$ with adjoint $\eta(X)^*=X^*\odot\id_F$.

\item\label{str3}
The map $X\odot Y\mapsto\eta(X)Y$ defines an isometry from the tensor product of $\sB^a(\cB,E)$ and $\sB^a(\cC,F)$ over $\sB^a(\cB)$ onto a strictly dense subset of $\sB^a(\cC,E\odot F)$.
\end{enumerate}
\eprop

\proof
For Part \ref{str1} see, for instance, \cite[Corollary 1.20]{MSS06}. Part \ref{str2} is general theory of tensor products.

For Part \ref{str3} let us choose a bounded approximate unit $\bfam{u_\lambda}_{\lambda\in\Lambda}$ for $\cB$. Then
\beqn{
\eta(X)Yc
~=~
\lim_\lambda\eta(X)u_\lambda Yc
~=~
\lim_\lambda(Xu_\lambda)\odot(Yc),
}\eeqn
where we made use of $u_\lambda y\to y$ in norm for all $y\in F$. (This follows from Part \ref{str1}, but may also easily be verified by three epsilons.) It follows that
\bmun{
\AB{\eta(X)Yc,\eta(X')Y'c'}
~=~
\lim_\lambda\AB{(Xu_\lambda)\odot(Yc),(X'u_\lambda)\odot(Y'c')}
~=~
\lim_\lambda\AB{Yc,u_\lambda^*\AB{X,X'}u_\lambda Y'c'}
\\
~=~
\lim_\lambda\AB{u_\lambda Yc,\AB{X,X'}u_\lambda Y'c'}
~=~
\AB{Yc,\AB{X,X'}Y'c'}
~=~
c^*\AB{Y,\AB{X,X'}Y'}c'.
}\emun
Clearly, when restricted to the subset $E\odot F$ of $\sB^a(\cB,E)\odot\sB^a(\cC,F)$, we obtain all maps of the form $c\mapsto cz$ for $z\in E\odot F$ that form a strictly dense subset of $\sB^a(\cC,E\odot F)$.\qed

\bdefi
By the \hl{strict tensor product} $\sB^a(\cB,E)\sodot\sB^a(\cC,F)$ we understand the space $\sB^a(\cC,E\odot F)$.
\edefi

The following corollary can be proved as Part \ref{str3}.

\bcor\label{spacecor}
For every correspondence $G$ from $\cC$ to $\cD$ (that may be viewed also as a correspondence $G$ from $\sB^a(\cC)$ to $\cD$ in a unique way) we have
\beqn{
\bfam{\sB^a(\cB,E)\sodot\sB^a(\cC,F)}\odot G
~=~
\sB^a(\cC,E\odot F)\odot G
~=~
E\odot F \odot G
}\eeqn
\ecor

\lf
Clearly, if $E$ is a correspondence from $\cA$ to $\cB$, then $\sB^a(\cC,E\odot F)$ is a correspondence from $\sB^a(\cA)$ to $\sB^a(\cC)$. In particular, if $E^\odot$ is a product system of correspondences over $\cB$, then the family of all $\sB^a(\cB,E_t)$ form a strict tensor product system of correspondences over $\sB^a(\cB)$. If this product system has a unital unit $\Xi^\odot$, then we may proceed as in Section \ref{psunisec}. So $\Xi_s\odot\id_{E_t}$ defines an inductive system of isometric embeddings $\sB^a(\cB,E_t)\rightarrow\sB^a(\cB,E_{s+t})$. From the inductive limit of this system we may extract a Hilbert \nbd{\cB}module
\beqn{
E_\infty
~=~
\bfam{\limind_{t\to\infty}\sB^a(\cB,E_t)}\odot\cB
}\eeqn
so that ~~$\limind_{t\to\infty}\sB^a(\cB,E_t)$ embeds as a strictly dense subset into $\sB^a(\cB,E_\infty)$. (Note that $\cB$ is a self-inverse Morita equivalence over $\cB$ and that the left action of $\cB$ extends to a strict left action of $\sB^a(\cB)$ on $\cB$ in the canonical way.) $\sB^a(\cB,E_\infty)$ fulfills
\beq{\label{strdil}
\sB^a(\cB,E_\infty)\sbar{\odot}\sB^a(\cB,E_t)
~=~
\sB^a(\cB,E_\infty)
}\eeq
and the usual associativity condition like \eqref{asscond}, so that $\vt_t(a)=a\odot\id_{E_t}$ defines an \nbd{E_0}semigroup on $\sB^a(\sB^a(\cB,E_\infty))$. But, we know that $\sB^a(\sB^a(\cB,E_\infty))$ is just $\sB^a(E_\infty)$. It is easy to show that that this \nbd{E_0}semigroup is strict and that its product system is nothing but $E^\odot$. The following proposition is slightly more general and implies what we just asserted in the special case $M=\cB$.

\bprop\label{strMoprop}
Let $M$ denote a Morita equivalence from $\cB$ to $\cC$ (so that $M$ carries a unique and strict extension of its left action to $\sB^a(\cB)$). Put $F_\infty:=E_\infty\odot M$. Then
\beqn{
F_\infty
~=~
F_\infty\odot(M^*\odot E_t\odot M)
}\eeqn
(via \eqref{strdil} and Corollary \ref{spacecor}) and $\theta_t(a):=a\odot\id_{M^*\odot E_t\odot M}$ defines a strict \nbd{E_0}semigroup on $\sB^a(F_\infty)$ whose product system is $F^\odot:=M^*\odot E^\odot\odot M$.
\eprop

\proof
The isomorphism $F_\infty=F_\infty\odot F_t$ is 
\bmun{
F_\infty
~=~
E_\infty\odot M
~=~
\sB^a(\cB,E_\infty)\odot M
~=~
\bfam{\sB^a(\cB,E_\infty)\sodot\sB^a(\cB,E_t)}\odot M
\\
~=~
E_\infty\odot E_t\odot M
~=~
E_\infty\odot M\odot M^*\odot E_t\odot M
~=~
F_\infty\odot F_t.
}\emun
The remaining statements follow as in the second half of the proof of Theorem \ref{psMethm} just the roles of $E^\odot$ and $F^\odot$ have now switched.\qed

\bthm\label{strunithm}
Let $E$ be a full correspondence over a unital \nbd{C^*}algebra $\cB$. Then there is a (necessarily full) Hilbert \nbd{\cB}module $F$ and a unital endomorphism of $\vt$ of $\sB^a(F)$ such that $F_\vt=E$.
\ethm

\proof
Denote by $E^\odot=\bfam{E^{\odot n}}_{n\in\N_0}$ the product system generated by $E$. We define $M_\infty(\cB)$ and $M_\infty(E_t)$ as the completions of the spaces of matrices with finitely many nonzero entries in the respective norm topologies and operations like in \eqref{matop}. To come to the setting of the preceding proposition we make up a dictionary.
\beqn{
\begin{array}{c|c}
\text{Propositions \ref{strMoprop}}&\text{here}
\\\hline
\cB		&M_\infty(\cB)
\\
\cC		&\cB
\\
M		&\cB^\infty
\\
E^\odot	&M_\infty(E^\odot):=\bfam{M_\infty(E_n)}_{n\in\N_0}
\\
F^\odot	&E^\odot
\\
F_\infty	&F
\\
\theta		&\vt
\end{array}
}\eeqn
In order to apply Proposition \ref{strMoprop} (providing us with the $F$ and the $\vt$ we seek according to the dictionary) it remains to show that $\sB^a(M_\infty(\cB),M_\infty(E))$ has a unit vector $\Xi$ (determining a unital unit $\Xi^\odot$ for the whole product system $M_\infty(E^\odot)$ as ingredient). But this can be done as in Proposition \ref{inftyprop} using, however, the ingredients from Lemma \ref{univeclem} (that is, $\el$ finite so that $\en=\#\N$ is sufficient) instead of those from Lemma \ref{Wuniveclem}.\qed

\bob\label{C*uniexob}
Also here the first part of Observation \ref{uniexob} remains true: $F\supset E^\infty$ contains a unit vector. The second half, $F$ can be chosen $\cB^\infty$, remains true at least if $E$ is countably generated. (This follows from the main result of Brown, Green and Rieffel \cite{BGR77}. We do not give any detail.)
\eob

\brem\label{nfrem}
Also here a correspondence $E$ over $\cB$ that should come from an endomorphism, necessarily must be also a correspondence over $\cB_E$. However, under the present assumptions we cannot simply replace $\cB$ with $\cB_E$ as in Remark \ref{nsfrem}, because $\cB_E$, in general, will be nonunital.

It is not difficult to write down endomorphisms or even continuous time \nbd{E_0}semigroups of $\sB^a(F)$ for a full Hilbert module $F$ over a nonunital \nbd{C^*}algebra. (Put $\cB:=C_0(0,\infty)$ and $F:=C_0((0,\infty),H)$ for some nonzero Hilbert space $H$. Then $F$ is a Hilbert \nbd{\cB}module with inner product $\AB{h,h'}(r):=\AB{h(r),h'(r)}$. Define $E_t=\cB$ as right Hilbert module, but with left action $b.x(r)=b(r+t)x(r)$. Then $F\odot E_t=F$ via $\SB{h\odot x}(r)=h(r+t)x(r)$ defines an \nbd{E_0}semigroup $\vt_t(a)=a\odot\id_t$ on $\sB^a(F)$ with product system $E_t$.) Nevertheless, without going into detail, we would like to emphasize that in many respects our motivation to study \nbd{E_0}semigroups on $\sB^a(F)$ via product systems (dilation theory!) lets appear as not very natural the case where $F$ is not full over a unital \nbd{C^*}algebra. Continuous product systems of correspondences over a unital \nbd{C^*}algebra always have unit vectors; see \cite[Lemma 3.2]{Ske07}.
\erem

\section{Representations}\label{repsec}

In this section we resolve Problem 2. We show that every faithful \nbd{C^*}correspondence admits a faithful nondegenerate representation (Theorem \ref{C*Hirthm}). This generalizes Hirshberg's result \cite{Hir05a} where the correspondence is required full. And we show that for a faithful \nbd{W^*}cor\-re\-spond\-ence the representation may be chosen normal (Theorem \ref{W*Hirthm}). Actually, we show first the result for the \nbd{W^*}case and, then, boil down the \nbd{C^*}case to the \nbd{W^*}case. The heart of the proof is the well-known statement that faithful representations of \nbd{W^*}algebras become unitarily equivalent when amplified suitably (Corollary \ref{ampisocor}). The main reason why we reproved that fact in Section \ref{W*uvsec} is to underline how simple a self-contained proof of Theorem \ref{W*Hirthm} actually is. The reduction of Theorem \ref{C*Hirthm} to Theorem \ref{W*Hirthm} remains somewhat tedious.

Let $G$ be a Hilbert space. A \hl{representation on $G$} of a correspondence $E$ over $\cB$ is a pair $(\pi,\eta)$ of maps $\pi\colon\cB\rightarrow\sB(G)$ and $\eta\colon E\rightarrow\sB(G)$ where $\pi$ is a representation of $\cB$ and $\eta$ is a bimodule map (that is, $\eta(bxb')=\pi(b)\eta(x)\pi(b')$) such that $\eta(x)^*\eta(y)=\pi(\AB{x,y})$. We always assume that $\pi$ is nondegenerate. The representation $(\pi,\eta)$ is \hl{nondegenerate} (or \hl{essential}), if also $\eta$ is nondegenerate, that is, if $\cls\eta(E)G=G$.

\brem\label{MSrem}
The nomenclature here differs, for instance, from Muhly and Solel \cite{MuSo98}, who call \it{covariant representation} a pair $(\pi,\eta)$ of completely contractive mappings fulfilling all conditions but $\eta(x)^*\eta(y)$ $=\pi(\AB{x,y})$. They call a covariant representation \it{isometric} if also $\eta(x)^*\eta(y)=\pi(\AB{x,y})$ holds, and they call an isometric covariant representation (that is, a representation in our sense) \it{fully coisometric} if $\cls\eta(E)G=G$.
\erem

We are done with Problem 2, if we can choose a faithful representation $\pi$ such that the induced representation $\rho^\pi_\cB$ of $\cB$ on $H:=E\odot G$ is unitarily equivalent to $\pi$, so that there exists a unitary $u\in\sB(G,H)$ such that $u\pi(b)=\rho^\pi(b)u$ for all $b\in\cB$. In that case, by setting $\eta(x)=u^*\eta^\pi(x)\in\sB(G)$ the pair $(\pi,\eta)$ is a faithful nondegenerate representation of $E$ on $G$. If, in the \nbd{W^*}case, $\pi$ is normal, then so is $\eta$ (see Section \ref{Irep}).

\bthm\label{W*Hirthm}
Every faithful \nbd{W^*}correspondence over a \nbd{W^*}algebra admits a normal faithful nondegenerate representation on a Hilbert space.
\ethm

\proof
Let $E$ be a \nbd{W^*}correspondence over a \nbd{W^*}algebra $\cB$. Choose a faithful normal nondegenerate representation $\pi\colon\cB\rightarrow\sB(G)$ of $\cB$ on a Hilbert space $G$. Then the induced representation $\rho^\pi_\cB$ on $H:=E\odot G$ is nondegenerate and normal. It is faithful because the left action of $\cB$ on $E$ is faithful. By Corollary \ref{ampisocor} there exists a Hilbert space $\eH$ such that the amplification $\pi\otimes\id_\eH$ of $\pi$ on $G\otimes\eH$ and the amplification $\rho^\pi_\cB\otimes\id_\eH$ of $\rho^\pi_\cB$ on $H\otimes\eH$ are unitarily equivalent. Obviously, $E\odot(G\otimes\eH)=H\otimes\eH$ so that $\rho^\pi_\cB\otimes\id_\eH$ is the representation $\rho^{\pi\otimes\id_\eH}_\cB$ of $\cB$ induced by $\pi\otimes\id_\eH$. By the discussion preceding the theorem we find a faithful normal nondegenerate representation $(\pi\otimes\id_\eH,\eta)$ of $E$.\qed

\bthm\label{C*Hirthm}
Every faithful \nbd{C^*}correspondence over a \nbd{C^*}algebra admits a faithful nondegenerate representation on a Hilbert space.
\ethm

\lf\noindent
\sc{Proof.}
Suppose $E$ is a correspondence over $\cB$ with a faithful left action. We are done, if we can choose a faithful (nondegenerate) representation $\pi\colon\cB\rightarrow\sB(G)$ in such a way that the induced representation $\rho^\pi_\cB$ on $H:=E\odot G$ extends to a normal and faithful representation of $\cB'':=\pi(\cB)''\subset\sB(G)$. This representation turns, then, $E'':=\ol{\eta^\pi(E)}^s$, the strong closure of the subset $\eta^\pi(E)$ in $\sB(G,H)$, into a von Neumann correspondence over $\cB''$ with faithful left action. (See Section \ref{comndrsec} for details.) We apply Theorem \ref{W*Hirthm} to $E''$ and obtain a (normal) faithful nondegenerate representation $\eta''$ of $E''$ on a Hilbert space. As $E$ is strongly dense in $E''$ (via $\eta^\pi$) and $\eta''$ is normal, also the restriction of $\eta:=\eta''\upharpoonright E$ to $E$ is nondegenerate.

Let ${\cB^*_1}^+:=\bCB{\vp\in\cB^*\colon\vp\ge0,\norm{\vp}\le1}$ and $E_1:=\bCB{x \in E\colon\norm{x}\le1}$. Suppose we can find a subset $S$ of ${\cB^*_1}^+$ that fulfills:
\begin{enumerate}
\item
For all $b\ne0$, there is a $\vp\in S$ such that $\vp(b^*b)\ne0$.

\item
For all $\vp\in S$ and $x\in E_1$, also $\vp\circ\AB{x,\bullet x}\in S$.

\item
For every $\vp\in S$, there exist $\psi\in S$ and $x\in E_1$ such that $\vp=\psi\circ\AB{x,\bullet x}$.
\end{enumerate}
We represent $\cB$ by $\pi=\bigoplus_{\vp\in S}\pi_\vp$ on $G=\bigoplus_{\vp\in S}G_\vp$ as the direct sum of all GNS-representations $(\pi_\vp,G_\vp\ni\gamma_\vp)$ to all elements $\vp$ in $S$. Then, by (1) this representation of $\cB$ is faithful. By (2) the induced representation of $\cB$  on $H:=E\odot G=\bigoplus_{\vp\in S}E\odot G_\vp$ extends to a normal representation of $\cB''\subset\sB(G)$. (Indeed, for $\vp\in S$, $x\in E_1$ and $\psi=\vp\circ\AB{x,\bullet x}\in S$ we observe that the subspace $H_{\vp,x}:=\cls\cB x\odot\gamma_\vp$ with the natural left action of $\cB$ is unitarily equivalent to the GNS-representation $\pi_\psi$ on $G_\psi$ by $v_{\vp,x}\colon bx\odot\gamma_\vp\mapsto b\gamma_\psi$. For $b''\in\cB''$ we simply define the action on an element $h\in H_{\vp,x}$ as $v_{\vp,x}^*b''v_{\vp,x}h$. It is easy to see that this extends as a well-defined representation of $\cB''$ on all of $H$ that is strongly continuous on bounded subsets and, therefore, normal.) And by (3) this representation of $\cB''$ is faithful.

For sequences $\bfam{\vp_n}_{n\in\N}$ in ${\cB^*_1}^+$ and $\bfam{x_n}_{n\in\N}$ in $E_1$ we denote
\beqn{
\vp^n
~:=~
\vp_n\circ\AB{x_1\odot\ldots\odot x_n,\bullet x_1\odot\ldots\odot x_n}.
}\eeqn
Recall that a \hl{subnet} of a net $\bfam{a_\lambda}_{\lambda\in\Lambda}$ is a net of the form $\bfam{a_{g(\mu)}}_{\mu\in M}$ for some \hl{cofinal} function $g\colon M\rightarrow\Lambda$ (that is, for every $\lambda\in\Lambda$ there is a $\mu_\lambda\in M$ such that $\mu\ge\mu_\lambda\Rightarrow g(\mu)\ge\lambda$). We define a suitable set $S$ by
\bmun{
S
~:=~
\Bset[\vp\in{\cB^*_1}^+]{\exists~\bfam{\vp_n}_{n\in\N}\subset{\cB^*_1}^+,\bfam{x_n}_{n\in\N}\subset E_1\text{~such that}
\\
\vp\text{~is the weak$^*$ limit of a subnet of~}\bfam{\vp^n}_{n\in\N}}.
}\emun

To show (1), let $b\in\cB$ with $\norm{b}=1$. Then choose $x_n\in E_1$ such that
\beqn{
\norm{bx_1\odot\ldots\odot x_n}^2
~\ge~
\frac{n+1}{2n}
}\eeqn
for all $n\in\N$, and choose states $\vp_n$ such that
\beqn{
\vp_n\bfam{\AB{bx_1\odot\ldots\odot x_n,bx_1\odot\ldots\odot x_n}}
~=~
\norm{bx_1\odot\ldots\odot x_n}^2.
}\eeqn
${\cB^*_1}^+$ is weak$^*$ compact so that the sequence $\bfam{\vp^n}_{n\in\N}$ has a weak$^*$ convergent subnet. Its limit $\vp$ is an element of $S$ that fulfills $\vp(b^*b)\ge\frac{1}{2}$.

For (2) and (3) let us fix an arbitrary element of $\vp\in S$ represented as weak$^*$ limit $\vp=\lim_\lambda\vp^{f(\lambda)}$ for some sequences $\bfam{\vp_n}_{n\in\N}$ in ${\cB^*_1}^+$ and $\bfam{x_n}_{n\in\N}$ in $E_1$, a directed set $\Lambda$ and a cofinal function $f\colon\Lambda\rightarrow\N$.

To show (2), choose $x\in E_1$. Then for $\psi_n=\vp_{n-1}, y_n=x_{n-1}$ $(n\ge2)$ and $\psi_1=0,y_1=x$ and the cofinal function $g(\lambda)=f(\lambda)+1$ we find that $\psi^{g(\lambda)}(b)=\psi_{g(\lambda)}\bfam{\AB{y_1\odot\ldots\odot y_{g(\lambda)},by_1\odot\ldots\odot y_{g(\lambda)}}}=\vp^{f(\lambda)}(\AB{x,bx})$ converges to $\vp(\AB{x,bx})$ for all $b\in\cB$. So, $\vp\circ\AB{x,\bullet x}=\lim_\lambda\psi^{g(\lambda)}\in S$.

To show (3), a candidate for $x$ is $x_1$. We put $y_n=x_{n+1},\psi_n=\vp_{n+1}$ $(n\in\N)$ and $g(\lambda)=\max(f(\lambda)-1,1)$. By weak$^*$ compactness, from the net $\bfam{\psi^{g(\lambda)}}_{\lambda\in\Lambda}$ we may choose a subnet $\bfam{\psi^{g\circ h(\mu)}}_{\mu\in M}$ converging weakly to a $\psi$. Clearly, the function $g\circ h\colon M\rightarrow\N$ is cofinal, so that $\psi\in S$. And $\psi$ fulfills $\psi\circ\AB{x_1,\bullet x_1}=\vp$, because the net $\bfam{\psi^{g\circ h(\mu)}\circ\AB{x_1,\bullet x_1}}_{\mu\in M}$ has the same limit as the subnet $\bfam{\vp^{f\circ h(\mu)}}_{\mu\in M}$ of $\bfam{\vp^n}_{n\in\N}$, namely, $\vp$. (These two nets are identical for all $\mu$ apart from those where $f(h(\mu))=1$.)\qedsymbol

\brem\label{nfaithrem}
If $E$ is a (\nbd{W^*})correspondence that admits a faithful (normal) nondegenerate representation $(\pi,\eta)$, then necessarily $E$ is faithful. (The induced representation $\rho^\pi_\cB$ is unitarily equivalent to the faithful representation $\pi$.) If $E$ is (strongly) full, then for $\pi$ being faithful it is sufficient (and necessary) that $\eta$ alone is faithful.

What happens, if we require only that $\eta$ is faithful in the case when $E$ is not necessarily (strongly) full? In this case, at least the restriction of $\pi$ to $\cB_E$ (respectively, $\ol{\cB_E}^s$) must be faithful and the left action of $\cB_E$ (respectively, $\ol{\cB_E}^s$) must be nondegenerate and faithful. It follows that $E$ must be faithful as (\nbd{W^*})correspondence over $\cB_E$ (respectively, $\ol{\cB_E}^s$). But in this case we can apply Theorem \ref{C*Hirthm} (\ref{W*Hirthm}) and obtain a faithful representation $(\eta,\pi)$ of the correspondence over the smaller algebra. The part $\pi$ of that representation may be extended to a (of course, in general not faithful) representation of $\cB$ and the extended $\pi$ together with $\eta$ gives rise to a (normal) nondegenerate not necessarily faithful representation where $\eta$ is faithful.
\erem

\bob\label{repampob}
The proof of Theorem \ref{W*Hirthm} shows that, if $\cB\subset\sB(G)$ is a von Neumann algebra, then the representation of a \nbd{W^*}correspondence $E$ over $\cB$ may be chosen to live on a multiple of the representation space $G$ with $\cB$ acting in the natural way as amplification. Also, in the \nbd{C^*}case the representation of $E$ will live on a suitable multiple of the representation space of the representation $\pi$ constructed on the proof of Theorem \ref{C*Hirthm}.
\eob

The following consequence of Theorems \ref{C*Hirthm}, \ref{W*Hirthm}, and Section \ref{psunisec} is a discrete time version of a result by Arveson and Kishimoto \cite{ArKi92} for \nbd{W^*}algebras. In Skeide \cite{Ske06p6} we prove a continuous time version for \nbd{C^*}modules. In Skeide \cite{Ske07p} we will use the same technique to give a completely different proof of \cite{ArKi92}.

\bthm\label{AKthm}
Every faithful strict (normal) unital endomorphism $\vt$ of $\sB^a(F)$ for some Hilbert module (\nbd{W^*}module) over $\cB$ is the restriction of an automorphism of some $\sB(H)$ containing $\sB^a(F)$ as (\nbd{W^*})subalgebra.
\ethm

\proof
We discuss only the (more difficult) \nbd{C^*}case. By making $\cB$ smaller, we assume that $F$ is full. Denote by $E=F^*\odot{_\vt}F$ the correspondence of $\vt$. Since $\vt$ is faithful, so is $E$. Applying Theorem \ref{C*Hirthm}, we obtain a faithful nondegenerate representation $(\pi,\eta)$ of $E$ on $G$. Define $H:=F\odot G$. Since $\pi$ is faithful, so is the embedding $a\mapsto a\odot\id_G$ from $\sB^a(F)$ into $\sB(H)$. Since $\eta$ is nondegenerate, the elements $\eta(y^*\odot_\vt z)g$ are total in $G$. By
\beqn{
x\odot\eta(y^*\odot_\vt z)g
~\longmapsto~
\vt(xy^*)z\odot g
}\eeqn
we define a unitary $u\in\sB(H)$. For every $a\in\sB^a(F)$, we find that
\bmun{
u(a\odot\id_G)(x\odot\eta(y^*\odot_\vt z)g)
~=~
u(ax\odot\eta(y^*\odot_\vt z)g)
\\
~=~
\vt(axy^*)z\odot g
~=~
(\vt(a)\odot\id_G)(\vt(xy^*)z\odot g)
~=~
(\vt(a)\odot\id_G)u(x\odot\eta(y^*\odot_\vt z)g).
}\emun
In other words, $u(a\odot\id_G)u^*=\vt(a)\odot\id_G$, so that the restriction of the inner automorphism $u\bullet u^*$ of $\sB(H)$ to the subalgebra $\sB^a(F)\odot\id_G\cong\sB^a(F)$ gives back $\vt$.\qed

\bcor
Every faithful (normal) nondegenerate endomorphism of a \nbd{C^*}~ (\nbd{W^*})algebra $\cB$ is the restriction of an inner automorphism of some $\sB(H)\supset\cB$ to $\cB$.
\ecor

\proof
If $\cB$ is a von Neumann algebra, apply Theorem \ref{AKthm} to $\sB^a(\cB)=\cB$. If $\cB$ is a \nbd{C^*}algebra, then the nondegenerate(!) homomorphism $\vt\colon\cB\rightarrow\cB\subset\sB^a(\cB)$ extends uniquely to a strict unital homomorphism of $\sB^a(\cB)$. Now we may apply Theorem \ref{AKthm}.\qed

\lf
As another application we prove that every normal faithful representation $(\sigma_0,\sigma)$ on $G$ of a faithful \nbd{W^*}correspondence $E$ over $\cB$ admits a \hl{nondegenerate extension} $(\tau_0,\tau)$ on $H\supset G$. (In the terminology of \cite{MuSo98}, every normal isometric covariant representation admits a unitary, that is isometric and fully coisometric, extension.) By this we mean that $(\tau_0,\tau)$ is normal nondegenerate (faithful) representation of $E$ such that $\tau(x)g=\sigma(x)g$ $\tau_0(b)g=\sigma_0(b)g$ for all $g\in G\subset H$.

Note that this is a stronger statement than existence of a nondegenerate dilation of $(\sigma_0,\sigma)$. Dilation would mean that the compression to $G$ gives back $(\sigma_0,\sigma)$. Existence of a nondegenerate dilation has been shown \cite{MuSo02}. As explained in \cite{Ske06p3}, existence also follows via the \it{commutant} (see Section \ref{comndrsec}) from the inductive limit construction described in Section \ref{psunisec} due to \cite{BhSk00,BBLS04}. It is known that the statement may fail for \nbd{C^*}correspondences, see Solel \cite[Example 5.16]{MuSo98}.

\bthm\label{ndthm}
Every normal faithful representation of a (faithful) \nbd{W^*}correspondence admits a nondegenerate extension.
\ethm

\proof
Let $E$ denote a faithful \nbd{W^*}correspondence over a \nbd{W^*}algebra $\cB$. Suppose $(\sigma_0,\sigma)$ is a normal isometric faithful covariant representation of $E$ on the Hilbert space $G$.

Then the Hilbert space $E\odot G$ is canonically isomorphic to the subspace $H:=\cls\sigma(E)G$ of $G$ and the induced representation $\rho^{\sigma_0}_\cB\colon b\mapsto b\odot\id_G$ on $E\odot G$ is unitarily equivalent to $\sigma_0\upharpoonright H$.

Since $E$ is faithful, by Theorem \ref{W*Hirthm} there exists a normal faithful nondegenerate representation $(\eta_0,\eta)$. By Observation \ref{repampob} the representation space may be chosen $G^{\#S}:=\bigoplus_{s\in S}G$ for some infinite set $S$ in such a way that $\eta_0=\id_\cB^{\#S}$.

Choose $s_0\in S$ and fix a bijection $\vp\colon S\rightarrow S\backslash\CB{s_0}$. For every Hilbert space $K$ we define a unitary $v_K\colon K\oplus K^{\#S}\rightarrow K^{\#S}$ by setting
\beqn{
v_K(k,\bfam{k_s}_{s\in S})
~=~
\bfam{k'_s}_{s\in S}
\text{~~~with~~~}
k'_{s_0}
~:=~
k,
~~~
k'_{\vp(s)}
~:=~
k_s
}\eeqn
(cf.\ Example \ref{SzNFex}). Denote by $H^\perp$ the orthogonal complement of $H$ in $G$. Of course, $G^{\#S}=H^{\#S}\oplus{H^\perp}^{\#S}$ in the obvious way. If we understand $\sigma(x)$ as an element in $\sB(G,H)$, then $(\eta_0,\tau)$ with
\beqn{
\tau(x)
~:=~
(v_H\oplus\id_{{H^\perp}^{\#S}})(\sigma(x)\oplus\eta^{\#S}(x))v_G^*
}\eeqn
defines a normal faithful nondegenerate representation of $E$ on $G^{\#S}$ that sends the subspace $G\cong G_{s_0}$ to $H_{s_0}\subset G_{s_0}$ and, on that subspace, gives back $\sigma$.\qed

\brem
The extension does not give an extension of the representation (in the sense of Definition \ref{psrepdef}) of the whole product system $\bfam{E^{\odot n}}_{n\in\N_0}$ generated by $E$. A ``semigroup'' version of this result has to wait for future investigation.
\erem

\section{Commutants: Endomorphisms \it{versus} representations}\label{comndrsec}

In this section we show that the \nbd{W^*}versions of our results, Theorem \ref{normunithm} and Theorem \ref{W*Hirthm}, are dual to each other in the sense of commutants of von Neumann correspondences. The commutant is a duality between a von Neumann correspondence over the von Neumann algebra $\cB\subset\sB(G)$ and its \it{commutant}, a von Neumann correspondence over the commutant $\cB'$ of $\cB$. In Theorem \ref{EE'switchthm} we will show that, under commutant, endomorphisms associated with a von Neumann correspondence $E$ are in correspondence (in a sense one-to-one) with representations of its commutant, $E''$. An endomorphism is unital, if and only if the corresponding representation is nondegenerate.

Von Neumann modules (Skeide \cite{Ske00b}) are the concrete operator analogues of \nbd{W^*}mod\-ules and Neumann correspondences (Skeide \cite{Ske03c,Ske06b}) are the concrete operator analogues of \nbd{W^*}correspondences. (As categories the two versions are equivalent.) Unlike the \nbd{W^*}version, for von Neumann modules there is a double commutant theorem and von Neumann correspondences posses a commutant. (The commutant was introduced in Skeide \cite{Ske03c}. Independently, Muhly and Solel \cite{MuSo04} have considered a \nbd{W^*}version, in which the \nbd{W^*}algebra, first, must be represented faithfully. In \cite{MuSo05} they generalized the construction to \nbd{\cA}\nbd{\cB}cor\-re\-spond\-ences.) We start by giving a very brief account on these subjects.

Let $\cB\subset\sB(G)$ be a von Neumann algebra acting nondegenerately on the Hilbert space $G$. Then every (pre-)Hilbert \nbd{\cB}module $E$ may be identified as a concrete operator \nbd{\cB}submodule of $\sB(G,H)$ (nondegenerate in the sense that $\cls EG=H$) via the representation $\eta:=\eta^{\id_\cB}$ from $G$ to $H$ induced by the identity representation $\id_\cB$ of $\cB$ on $G$; see Section \ref{Irep}. Following Skeide \cite{Ske00b}, we say $E$ is a \hl{von Neumann \nbd{\cB}module}, if $E$ is strongly closed in $\sB(G,H)$.

One may show that $E$ is a von Neumann module, if and only if $E$ is self-dual, that is, if and only if $E$ is a \nbd{W^*}module; see \cite{Ske00b,Ske05c}. For a fixed von Neumann algebra $\cB$ the category of von Neumann \nbd{\cB}modules and the category of \nbd{W^*}modules over $\cB$ are, therefore, equivalent. (The morphisms are, in both cases, the adjointable maps.) Fixing an equivalence between the category of \nbd{W^*}algebras and the category of von Neumann algebras, also the category of von Neumann modules and the category of \nbd{W^*}modules are equivalent. (The morphisms are the ternary morphisms; see Abbaspour and Skeide \cite{AbSk07} for details.)

\brem\label{dscrem}
The point about von Neumann modules is that it is easier to obtain them (from pre-Hilbert modules over a von Neumann algebra) than \nbd{W^*}modules. Simply take strong closure. In the sequel, we will learn another possibility that is completely algebraic and parallels the operation of taking the double commutant of an operator \nbd{*}algebra in order to obtain a von Neumann algebra; see Remark \ref{dctrem}.
\erem

We identify $\sB^a(E)$ as a subalgebra of $\sB(H)$ via the induced representation $\rho^{\id_\cB}$. Clearly, if $E$ is a von Neumann module, then $\sB^a(E)$ is a von Neumann subalgebra of $\sB(H)$. When $E$ is also a correspondence over $\cB$ such that the canonical representation $\rho\colon\cB\rightarrow\sB^a(E)\rightarrow\sB(H)$ is normal, then we say $E$ is a \hl{von Neumann correspondence}. (So $E$ is a von Neumann correspondence, if and only if it is also a \nbd{W^*}correspondence. Once more there are equivalences of von Neumann categories and \nbd{W^*}categories, with and without fixing the algebra in question.) We refer to $\rho$ as the \hl{Stinespring representation} of $\cB$.

On $H$ there is a second (normal nondegenerate) representation, namely, the so-called \hl{commutant lifting} $\rho'$ of $\cB'$ defined as $\rho'(b')=\id_E\odot b'$. It is not difficult to show that the intertwiner space $C_{\cB'}(\sB(G,H)):=\CB{x\in\sB(G,H)\colon \rho'(b')x=xb'~(b'\in\cB')}$ is a von Neumann \nbd{\cB}module (see \cite{Rie74a}) and that $E$ is a von Neumann module, if and only if $E=C_{\cB'}(\sB(G,H))$ (see \cite{Ske05c}). Less obvious is the converse statement: If $\rho'$ is a normal nondegenerate representation of $\cB'$ on a Hilbert space $H$, then the von Neumann \nbd{\cB}module $E:=C_{\cB'}(\sB(G,H))$ acts nondegenerately on $G$ (see \cite[Lemma 2.10]{MuSo02}), that is, $E\odot G=H$ via $x\odot g=xg$. Clearly, the commutant lifting for that $E$ is the $\rho'$ we started with. The fact that the correspondence between von \nbd{W^*}modules over $\cB$ and representations of $\cB'$ (in standard representation) is an equivalence of categories, has been observed in Baillet, Denizeau and Havet \cite{BDH88}. (\cite[Theorem 2.2]{BDH88} is, actually, between \nbd{W^*}correspondences and correspondences in the sense of Connes \cite{Con80p}. One has to put the algebra acting from the left to $\C$.) A version as a bijective functor (between \it{concrete} von Neumann \nbd{\cB}modules and representations of $\cB'$) is due to Skeide \cite{Ske06b}.

\brem\label{dctrem}
If $E$ is only a pre-Hilbert module over the von Neumann algebra $\cB$, then $\ol{E}^s$ is just $C_{\cB'}(\sB(G,H))$ and provides us with the minimal self-dual extension of $E$ in the sense of Paschke \cite{Pas73}; see \cite{Rie74a,Ske05c}. This is the \it{double commutant theorem} for von Neumann modules.
\erem

A von Neumann \nbd{\cB}module is strongly full, if and only if the commutant lifting $\rho'$ is faithful. 

\bcor\label{freecor}
If $F$ is a strongly full von Neumann module, then there is a cardinal number $\en$ such that $\ol{F^\en}^s\cong\ol{\cB^\en}^s$. In particular, if $\vt$ is a unital normal endomorphism of $\sB^a(F)$ with associated von Neumann correspondence $E$, then the amplification gives a unital normal endomorphism $\vt^\en$ on $\sB^a(\ol{F^\en}^s)=\sB^a(\ol{\cB^\en}^s)=\cB\otimes\sB(\C^\en)$, whose associated von Neumann correspondence is $E$, too.
\ecor

\proof
The first statement is a simple consequence of Corollary \ref{ampisocor} and the observation that the correspondence between von Neumann modules and their commutant liftings respects direct sums (of arbitrary cardinality). The second statement follows from the (easy to proof) fact that \nbd{E_0}semigroups with the same associated product system may be added.\qed

\lf
As a curiosity we reprove a well-known result (see \cite[Theorem 8.15]{Rie74a} and its footnote) about when two \nbd{W^*}algebras are Morita equivalent.

\bcor\label{Mecor}
Let $\cA$ and $\cB$ denote two \nbd{W^*}algebras. Then the following conditions are equivalent:
\begin{enumerate}
\item
$\cA$ and $\cB$ are Morita equivalent.

\item
$\cA$ and $\cB$ admit faithful normal nondegenerate representations ~$\rho\colon\cA\rightarrow\sB(H)$~ and ~$\pi\colon\cB\rightarrow\sB(G)$~ with isomorphic commutants $\rho(\cA)'\cong\pi(\cB)'$.

\item
There is a Hilbert space $\eH$ such that $\cA\otimes\sB(\eH)$ and $\cB\otimes\sB(\eH)$ are isomorphic.

\end{enumerate}
\ecor

\proof
$1$ $\Rightarrow$ $2$.
Suppose $E$ is a Morita equivalence from $\cA$ to $\cB$. Choose a faithful normal nondegenerate representation $\pi$ of $\cB$ on $G$. Put $\cB':=\pi(\cB)'$ and define as usual the commutant lifting $\rho'$ of $\cB'$ on $H:E\odot G$. Since $E$ is a Morita equivalence, $\rho'$ is faithful and $\rho^\pi$ is an isomorphism onto $\sB^a(E)=\rho'(\cB')'\subset\sB(H)$. In other words, $\rho(\cA)'=\rho'(\cB')\cong\cB'=\pi(\cB)'$.

$2$ $\Rightarrow$ $3$.
Suppose we have two representations $\pi$ and $\rho$ as stated. Then Corollary \ref{ampisocor} provides us with a Hilbert space $\eH$ such that $\rho(\cA)'\otimes\id_\eH$ and $\pi(\cB)'\otimes\id_\eH$ are unitarily equivalent. Thus, they have isomorphic commutants $\rho(\cA)\otimes\sB(\eH)$ and $\pi(\cB)\otimes\sB(\eH)$ so that also $\cA\otimes\sB(\eH)$ and $\cB\otimes\sB(\eH)$ are isomorphic.

$3$ $\Rightarrow$ $1$.
~$\ol{\cA^{\idim\eH}}^s$ is a Morita equivalence from $\cA\otimes\sB(\eH)$ to $\cA$ and ~$\ol{\cB^{\idim\eH}}^s$ is a Morita equivalence from $\cB\otimes\sB(\eH)$ to $\cB$. If $\cA\otimes\sB(\eH)$ and $\cB\otimes\sB(\eH)$ are isomorphic, then the tensor product ~${\ol{\cA^{\idim\eH}}^s}^*\sodots\ol{\cB^{\idim\eH}}^s$ over $\cA\otimes\sB(\eH)\cong\cB\otimes\sB(\eH)$ makes sense and is a Morita equivalence from $\cA$ to $\cB$.\qed

\lf
Summarizing, we have a one-to-one correspondence between von Neumann \nbd{\cB}modules and representations of $\cB'$ and a one-to-one correspondence between von Neumann correspondences $E$ over $\cB$ and pairs of representations $(\rho,\rho',H)$ of $\cB$ and $\cB'$ with mutually commuting range. In the latter picture of correspondences as two representations nobody prevents us from exchanging the roles of $\cB$ and $\cB'$. In that way, we obtain a further von Neumann correspondence, namely
\beqn{
E'
~:=~
C_\cB(\sB(G,H))
~:=~
\bCB{x'\in\sB(G,H)\colon\rho(b)x'=x'b~(b\in\cB)},
}\eeqn
this time over $\cB'$ with left action of $\cB'$ via $\rho'$. This duality between $E$ and its \hl{commutant} $E'$ was mentioned in \cite{Ske03c}. See Skeide \cite{Ske06b} for definitions (\it{concrete} von Neumann correspondences) where the commutant becomes, really, a bijective functor.

We are now in a position to formulate the theorem about the relation between Problem 1 and Problem 2 for von Neumann correspondences. But first let us recall that a von Neumann correspondence $E$ over $\cB$ is strongly full, if and only if the left action of $\cB'$ on the commutant $E'$ defines a faithful representation of $\cB'$ on $E'$. If $(\pi',\eta')$ is a faithful nondegenerate representation of $E'$, so that the left action of $\cB'$ on $E'$ is faithful, then necessarily $E$ is strongly full.

\bthm\label{EE'switchthm}
Let $E$ be a von Neumann correspondence over a von Neumann algebra $\cB\subset\sB(G)$ and $E'$ its commutant. Then the following conditions are equivalent.
\begin{enumerate}
\item
$E$ is the correspondence of a normal unital endomorphism $\vt$ of $\sB^a(F)$ for some strongly full von Neumann \nbd{\cB}module $F$.

\item
$E'$ admits a faithful normal nondegenerate representation $(\pi',\eta')$ on a Hilbert space $K$.
\end{enumerate}
Moreover, if either of the conditions is fulfilled, then $E$ is strongly full or, equivalently, the left action of $\cB'$ on $E'$ is faithful.
\ethm

\proof
Let $(\rho,\rho',H)$ be the triple that determines $E$ as $C_{\cB'}(\sB(G,H))$ and $E'$ as $C_\cB(\sB(G,H))$.

Suppose that $F$ is a strongly full von Neumann \nbd{\cB}module and that $\vt$ is a normal unital endomorphism of $\sB^a(F)$ such that $F=F\sodots E$ and $\vt(a)=a\odot\id_E$. (As $F$ is strongly full, $E$ is uniquely determined by these properties and necessarily $E$ is itself strongly full.) Put $K:=F\odot G$. As $F=F\sodots E$ we have $K=F\odot E\odot G$. (If the last factor in a tensor product is a Hilbert space, then norm closure is sufficient.) By construction we have $E\odot G=\cls EG=H=\cls E'G=E'\odot G$. (Note that $\cls EG=H=\cls E'G$ is true equality of Hilbert spaces. The equalities $E\odot G=\cls EG$ and $E'\odot G=\cls E'G$ are by canonical isomorphism.) We find
\beqn{
F\odot G
~=~
K
~=~
F\odot E'\odot G.
}\eeqn
There are several ways to understand why $\eta'(x')\colon y\odot g\mapsto y\odot x'\odot g$ is a well-defined element of $\sB(K)$. One is that $\eta'(x')=\id_F\odot x'$ where $x'$ is considered a \nbd{\cB}\nbd{\C}linear operator from $G$ to $H=E'\odot G$. Let $\pi'$ denote the (normal!) commutant lifting of $\cB'$ on $K=F\odot G$. We leave it as an instructive exercise to check that $(\pi',\eta')$ is a representation of $E'$ on $K$. Obviously this representation is nondegenerate. It is normal, because $\pi'$ is normal. It is faithful because $F$ is strongly full.

Suppose now that $(\pi',\eta')$ is a faithful normal nondegenerate representation of $E'$ on $K$. We put $F:=C_{\cB'}(\sB(G,K))$. As $\pi'$ is faithful, $F$ is strongly full. Again
\beqn{
F\odot G
~=~
K
~=~
F\odot E'\odot G
}\eeqn
now via $\eta'(x')(y\odot g)\mapsto y\odot x'\odot g$. (Note that the set $\eta'(E')F\odot G$ is total in $K$, because $\eta'$ is nondegenerate.) Again we substitute $E'\odot G=H=E\odot G$ so that $F\odot G=F\odot E\odot G$. The action of $b'\in\cB'$ on these spaces is the same. To see this we observe, first, that $b'(y\odot x\odot g)=y\odot x\odot b'g=y\odot\rho'(b')(x\odot g)$. Then, writing a typical element of $H=E\odot G$ not as elementary tensor $x\odot g$ but as elementary tensor $x'\odot g$ and recalling that the action of $b'$ on $x'\odot g$ is just $\rho'(b')$, we find
\beqn{
b'(y\odot x'\odot g)
~=~
y\odot\rho'(b')(x'\odot g)
~=~
y\odot b'x'\odot g
~=~
\eta'(b'x')(y\odot g)
~=~
\pi'(b')\eta'(x')(y\odot g).
}\eeqn
As the commutant liftings on $F\odot G$ and on $F\odot E\odot G$ coincide, also the modules $F$ and $F\sodots E$ (being intertwiner spaces for the same commutant lifting) must coincide and $\vt(a)=a\odot\id_E$ induces a unital normal endomorphism of $\sB^a(E)$. Once again, as $F$ is strongly full, a correspondence $E$ is determined uniquely by these properties, so that $F^*\sodots{_\vt}F$ gives us back $E$.\qed

\brem
Muhly and Solel \cite{MuSo99} have constructed from a nondegenerate representation $(\pi',\eta')$ on $K$ an endomorphism of $\pi'(\cB')'\subset\sB(K)$. Taking into account that this algebra coincides exactly with our $\sB^a(F)\subset\sB(K)$, puts into perspective the second part of the proof of Theorem \ref{EE'switchthm} with the result from \cite{MuSo99}. In fact, the constructions of the endomorphism are very much the same, except that we have added the construction of $F$ and the interpretation of the algebra on which the endomorphism acts as $\sB^a(F)$. This considerably facilitates understanding why everything is well-defined.
\erem

\bex
Suppose $E=H$ is a Hilbert space of dimension $n=2,3,\ldots,\infty$. Then the commutant $H'$ of $H$ is isomorphic to $H$ and we recover the well-known fact that representations of the Cuntz algebra $\cO_n$ correspond to endomorphisms of index $n$ of $\sB(K)$, and that nondegenerate representations correspond to unital endomorphisms. Note that the isomorphism $H\cong H'$ is by no means a trivial issue. One may see this by looking at the discrete product systems generated by $H$ and $H'$, respectively. One is the commutant of the other, but their product system structures are anti-isomorphic. This is the same relation as that between the Bhat system and the Arveson system constructed from an \nbd{E_0}semigroup on $\sB(K)$; see Skeide \cite{Ske07p1}.
\eex

We give now a version of Theorem \ref{EE'switchthm} for a whole product system. The following definition (from \cite{MuSo04}, but in a different terminology; see Remark \ref{MSrem}) extends suitably the definition of a representation of a single correspondence to the definition of a representation of a whole product system.

\bdefi\label{psrepdef}
A \hl{representation} of a product system $E^\odot$ of correspondences over a \nbd{C^*}al\-ge\-bra $\cB$ is a pair $(\pi,\eta)$ where $\pi$ is a nondegenerate representation of $\cB$ on a Hilbert space $K$ and $\eta=\bfam{\eta_t}_{t\in\bS}$ is a family such that each $(\pi,\eta_t)$ is a representation of $E_t$ on $K$ and such that
\beq{\label{psrep}
\eta_{s+t}(x_s\odot y_t)
~=~
\eta_s(x_s)\eta_t(y_t).
}\eeq
A representation is \hl{nondegenerate}, if every $(\pi,\eta_t)$ is nondegenerate. In case of product systems of \nbd{W^*}correspondences we require that $\pi$ (and, therefore, every $(\pi,\eta_t)$) is normal.
\edefi

Suppose $\eta=\bfam{\eta_t}_{t\in\bS}$ is a family of mappings fulfilling \eqref{psrep} and the isometricity condition $\eta_t(x_t)^*\eta(y_t)=\eta_0(\AB{x_t,y_t})$. It is easy to see that $(\eta_0,\eta)$ is a representation.

Speaking about a whole product system instead of a single correspondence, Theorem \ref{EE'switchthm} remains true (with practically no changes in the proof, apart from a view more indices) for product systems of von Neumann correspondences indexed by $\N_0$ or $\R_+$. We phrase it here.

\bthm\label{PSswitchthm}
Let $E^{\sodots}$ be a product system of von Neumann correspondences over a von Neumann algebra $\cB\subset\sB(G)$. Then also the commutant $E'^{\sodots}=\bfam{E'_t}_{t\in\bS}$ possesses a canonical structure of a product system. Suppose that all $E_t$ are strongly full or, equivalently, that all $E'_t$ are faithful. Then there is a one-to-one correspondence between normal \nbd{E_0}semigroups $\vt$ associated with $E^{\sodots}$ (acting on the operators of a necessarily strongly full von Neumann \nbd{\cB}module) and nondegenerate normal faithful representations $(\pi',\eta')$ of $E'^{\sodots}$.
\ethm

\proof
Just do for every couple $\vt_t$ and $\eta_t$ what we did in the proof of Theorem \ref{EE'switchthm} for single mappings, and verify the additional conditions. This proceeding also reveals automatically how the product system structure of the commutant of a product system must be defined.\qed

\brem
The theorem has two extensions. The first is to the nonfull case. Here, by Remark \ref{nsfrem}, we must require that $\ol{\cB_{E_t}}^s$ is stationary for $t>0$ and acts nondegenerately on all $E_t$. (Recall from Section \ref{psunisec} that $E_0:=\cB$ is defined by hand.) We may phrase an equivalent condition on the $E'_t$, following Remark \ref{nfaithrem}. Dropping strong fullness, on the commutant side this leads to possibly non faithful $\pi'$ where, however, still every $\eta'_t$ is injective. All this can be proved very simply, by restricting $\cB$ to the smaller algebra $\ol{\cB_{E_t}}^s$ (acting nondegenerately on the subspace $\cls\cB_{E_t}G$ of $G$) and its commutant. Then we are in the strongly full case.

The second extension is to \hl{\nbd{E}semigroups}, that is, to semigroups of not necessarily unital endomorphisms. (The definition of the product system associated with an \nbd{E}semigroup on $\sB^a(E)$ is the same. The only difference is that now we do no longer obtain an isomorphism $E\odot E_t\rightarrow E$ but only an isometry onto the subspace $\vt_t(\U)E$ of $E$.) On the commutant side this leads to possibly degenerate representations. In this setting we are no longer sure that $\ol{\cB_{E_t}}^s$ is stationary for $t>0$, so we possibly leave also the strongly full case. This time $\eta'_t$ need no longer be injective. Anyway, also in this case we remain with a one-to-one correspondence of \nbd{E}semigroups associated with $E^{\sodots}$ and normal representations of $E'^{\sodots}$.
\erem

\section{Examples}\label{exsec}

In this section we discuss for two examples what Theorem \ref{normunithm} asserts. The first example discusses the correspondence in Example \ref{non1ex}. The reader might object that this correspondence is a Morita equivalence and that, therefore, the endomorphism granted by Theorem \ref{normunithm} is an automorphism. However, this is the simplest nontrivial example possible, and the discussion is already quite involved. The second example is a correspondence of a proper endomorphism. In the end of each example we discuss (due to space reasons only very briefly) the meanings of Theorems \ref{W*Hirthm} and \ref{EE'switchthm}.

\bex
As in Example \ref{non1ex} we put $\cB=\tMatrix{\C&0\\0&M_2}\subset M_3$ and $E=\tMatrix{0&\C_2\\\C^2&0}\subset M_3$. The operations of the correspondence $E$ over $\cB$ are those inherited from $M_3$. This remains even true for the tensor product:
\beqn{
x\odot y
~=~
xy
~\in~
E\odot E
~=~
\cB.
}\eeqn
In particular,
\beqn{
E_n
~:=~
E^{\odot n}
~=~
\begin{cases}
\cB&\text{$n$ even,}
\\
E&\text{$n$ odd.}
\end{cases}
}\eeqn
Fortunately, the structure of Hilbert \nbd{\cB}modules $F$ is not much more complicated than that of Hilbert spaces and we still can say in advance how automorphisms of $\sB^a(F)$ may look like. In particular, we can say when an automorphism is associated with the correspondence $E$.

Let $p_1=\tMatrix{1&0\\0&0}$ and $p_2=\tMatrix{0&0\\0&\U}$ denote the two nontrivial central projections in $\cB$. Every Hilbert \nbd{\cB}module $F$ decomposes into the direct sum $F=F_1\oplus F_2$ with $F_i=Fp_i$. The summand $F_1$ has inner product in $\tMatrix{\C&0\\0&0}$. We may identify it with a Hilbert space $\eH_1$. The summand $F_2$ has inner product in $\tMatrix{0&0\\0&M_2}$. Its structure is therefore that of  a Hilbert \nbd{M_2}module. A short computation shows that
\beqn{
F_2
~=~
F_2\odot M_2
~=~
F_2\odot\C^2\odot\C_2
~=~
\eH_2\odot\C_2
~=~
\eH_2\otimes\C_2,
}\eeqn
where we defined the Hilbert space $\eH_2:=F_2\odot\C^2$ and where we used in the last step that there is no difference between the interior tensor product $\odot$ over $\C$ and the exterior tensor product $\otimes$.

We note that $F$ is also a \nbd{W^*}module. Also most tensor products we write down in the sequel are strongly closed if they are norm closed.

An operator $a$ on $F$ cannot mix the components in $F_1$ and in $F_2$. (To see this simply multiply with $p_i$ from the right and use right linearity of $a$.) Therefore, $a$ decomposes as $a=a_1\oplus a_2$ where each $a_i$ is an operator on $F_i$ alone. $a_1$ can be any element in $\sB(\eH_1)$, while $a_2$ must be an element in $\sB(\eH_2)$ that acts on $F_2=\eH_2\otimes\C_2$ as $a_2\otimes\id_{\C_2}$. (To see the latter statement we may, for instance, observe that tensoring with $\C_2$ is an operation of Morita equivalence so that $F_2$ and $\eH_2$, indeed, have the same operators.) We find $\sB^a(F)=\sB(\eH_1)\oplus\sB(\eH_2)$.

It is easy to check that an automorphism of $\sB^a(F)$ either sends  $\sB(\eH_i)$ onto $\sB(\eH_i)$ or sends $\sB(\eH_1)$ onto $\sB(\eH_2)$ and \it{vice versa}. The first type is simply implemented by two unitaries $u_i\in\sB(\eH_i)$. It is, therefore, conjugate to the identity automorphism and the associated correspondence is $\cB$.  In order to have the second case necessarily $\eH_1$ and $\eH_2$ are isomorphic, to a Hilbert space $\eH$ say, and the action of the automorphism is exchange of the two copies of $\sB(\eH)$ plus, possibly, an automorphism of the first type. This second case is, thus, simply the flip $\f(a_1\oplus a_2)=a_2\oplus a_1$ on $\sB(\eH)\oplus\sB(\eH)$ (up to conjugation with a unitary in $\sB(\eH)\oplus\sB(\eH)$).

We claim that the correspondence associated with the flip is $E$. We show this by giving an isomorphism from $F\odot E$ to $F$ that implements the flip as $a\mapsto a\odot\id_E$ and appeal to the uniqueness of the correspondence inducing $\f$. Indeed, one checks easily that
\beqn{
\sMatrix{h_1\\h_2\otimes v^*}\odot\sMatrix{0&v_2^*\\v_1&0}
~\longmapsto~
\sMatrix{h_2\AB{v,v_1}\\h_1\otimes v_2^*}
~~~~~~~~~~~~
(h_1,h_2\in\eH;v,v_1,v_2\in\C^2)
}\eeqn
defines a surjective isometry. Moreover, choosing an arbitrary unit vector $e\in\C^2$ we see that $(a_1\oplus a_2)\odot\id_E$ acting on
\beqn{
\sMatrix{h_1\\h_2\otimes v^*}
~=~
\sMatrix{h_2\\h_1\otimes e^*}\odot\sMatrix{0&v^*\\e&0}
}\eeqn
gives
\beqn{
\sMatrix{a_1h_2\\a_2h_1\otimes e^*}\odot\sMatrix{0&v^*\\e&0}
~=~
\sMatrix{a_2h_1\\a_1h_2\otimes v^*}
~=~
\f(a_1\oplus a_2)\sMatrix{h_1\\h_2\otimes v^*}.
}\eeqn

The discussion shows that a Hilbert \nbd{\cB}module $F$ with an endomorphism on $\sB^a(F)$ that has $E$ as associated correspondence must have the form $F=\eH\oplus(\eH\otimes\C_2)$ and that the endomorphism is the flip $\f$ on $\sB^a(F)=\sB(\eH)\oplus\sB(\eH)$ up to unitary equivalence in $\sB^a(F)$. That is, the possible endomorphisms associated with $E$ are simply classified by the dimension of $\eH$.

We ask now which of them can be obtained by the steps used in the proof of Theorem \ref{normunithm}. The answer is simple: $E$ does not have unit vectors, but $E^2$ has. As the cardinality that occurs in Lemma \ref{Wuniveclem} is $\el=2$, the minimal cardinality $\en$ in Proposition \ref{inftyprop} is simply countably infinite, which we denote $\en=\infty$. A unit vector $\Xi\in M_\en(E)$ gives rise to an isometry $\Xi\odot\Xi\in M_\en(E)\odot M_\en(E)=M_\en(\cB)$ that must be proper. Example \ref{SzNFex} tells us that inductive limit over the even half $M_\en(E)^{\odot 2n}$ will be an infinite-dimensional space. Therefore, $\eH$ cannot be finite-dimensional. It will simply have $\dim\eH=\en$. For $\en=\infty$ it is separable, otherwise it is nonseparable.

Let us now calculate the commutant of $E$. To that goal we consider $\cB\subset M_3=\sB(\C^3)$ as von Neumann algebra acting on $G=\C^3=\tMatrix{\C\\\C^2}$. We observe that $E\odot G=\tMatrix{0&\C_2\\\C^2&0}\odot\tMatrix{\C\\\C^2}=\tMatrix{\C\\\C^2}=G$ via $\tMatrix{0&x^*\\y&0}\odot\tMatrix{\lambda\\z}=\tMatrix{\AB{x,z}\\y\lambda}$. The Stinespring representation is just the identity representation. We find $E'=\cB'=\tMatrix{\C&0\\0&\C\U}$ and the identity $E'\rightarrow\sB(\C^3)$ is a normal faithful nondegenerate representation. The same is true for the identity representation of $E$ itself. So, as far as representations are concerned neither $E$ nor its commutant $E'$ yield interesting results. The only approximately noteworthy fact is that the commutant lifting for $E$ is the flip $\tMatrix{\lambda&0\\0&\mu\U}\mapsto\tMatrix{\mu&0\\0&\lambda\U}$.
\eex

\bex\label{lastEX}
We give now an example of a correspondence without unit vector, that comes from a proper endomorphism. Moreover, no tensor power of this correspondence admits a unit vector.

We consider the von Neumann algebra $\cB=\ol{\bigoplus_{n\in\N}M_n}^s$ acting on $G=\bigoplus_{n\in\N}\C^n$. Recall that $M_{nm}=\C^n\otimes\C_m$ is a von Neumann correspondence from $M_n$ to $M_m$ (actually, a Morita equivalence) that may also be considered as a correspondence over $\cB$. As $E$ we choose the von Neumann \nbd{\cB}correspondence direct sum
\beqn{
E
~:=~
\C\oplus\soplus_{n\in\N}\C^n\otimes\C_{n+1}.
}\eeqn
Here $\cB$ acts on direct summands of $E$ from either side with that direct summand $M_n$ that fits the correct dimension. That is, $M_1$ acts from the left on the summands $\C$ and $\C^1\otimes\C_2=\C_2$ but from the right only on $\C$. It is easy to check that
\beqn{
E^{\sodots m}
~:=~
\C\oplus\C_2\oplus\ldots\oplus\C_m\oplus\soplus_{n\in\N}\C^n\otimes\C_{n+m}.
}\eeqn
All $E^{\sodots m}$ are strongly full but none of them has a unit vector.

$E$ is not a Morita equivalence, so it must come from a proper endomorphism. To understand which endomorphisms could be associated with $E$, we analyze the general structure of a von Neumann \nbd{\cB}module $F$ and look for which (strongly full) $F$ we can write down an isomorphism $F\sodots E=F$. According to the minimal ideals $M_n$ in $\cB$, also $F$ decomposes into a direct sum of von Neumann \nbd{M_n}modules $F_n$. Every $F_n$ must have the form $\eH_n\otimes\C_n$ for some Hilbert space $\eH_n$. Of course, $\sB^a(F)=\ol{\bigoplus_{n\in\N}\sB(\eH_n)}^s$. A short computation yields that
\beqn{
F\sodots E
~:=~
\eH_1\oplus\soplus_{n\in\N}\eH_n\otimes\C_{n+1}.
}\eeqn
Therefore, $F\sodots E\cong F$, if and only if $\eH_n=\eH$ for all $n\in\N$. Another computation shows that the endomorphism induced by this isomorphism acts on $\sB^a(F)$ as $\vt(a_1,a_2,\ldots)=(a_1,a_1,a_2,\ldots)$. It is nothing but the unitalization of the one-sided shift on $\sB^a(F)$. As our construction of the inductive limit runs through a countable inductive system of proper isometries, an $\eH$ coming from our construction must be infinite-dimensional and separable. Note that in this case $F$ has a unit vector, while if $\eH$ is finite-dimensional, then $F$ fails to have a unit vector.

Clearly, the commutant of $\cB$ is $\cB'=\soplus_{n\in\N}\C\U_n$. Denote by $\s(z_1,z_2,z_3,\ldots):=(z_2,z_3,\ldots)$ the left shift on $\cB'$. We invite the reader to check that the commutant of $E$ is
\beqn{
E'
~=~
M_1\oplus{_{\text{\tiny$\sS$}}}\cB'
}\eeqn
and that the maps of the representation $(\pi',\eta')$ of $E'$ on $F\odot G=\bigoplus_{n\in\N}\eH$ granted by Theorem \ref{EE'switchthm} simply let act $(z_1,z_2,\ldots)\in\cB'$ and $\lambda\oplus(z_1,z_2,\ldots)=(\lambda,z_1,z_2,\ldots)\in E'$ component-wise on $\bigoplus_{n\in\N}\eH$.

$E$ is also a faithful von Neumann correspondence. The steps in the proof of Theorem \ref{W*Hirthm} to be carried out explicitly are very plain. Indeed, $E\odot G=\C\oplus G$ (with $\cB$ acting with its \nbd{1}component on the summand $\C$). Clearly, the components in the multiple $(\C\oplus G)^\infty$ can be rearranged easily to give a unitary equivalence with $G^\infty$ (including the respective actions of $\cB$ on these spaces). The von Neumann \nbd{\cB'}module induced by the representation of $\cB$ on $G^\infty$ is simply $F'=(\cB'^\infty)''$. The identification $F'\sodots E'=F'$ granted by the theorem simply identifies the infinitely many components $\C$ and $\C^1\subset\cB'$ contained $F'\sodots E'=\ol{(\C\oplus\cB')^\infty}^s$ with the infinitely many components of $\C^1\subset\cB'$ contained in $F'=(\cB'^\infty)''$. The remaining components of $\cB'$ remain untouched. The endomorphism $\vt'$ on $\sB^a(F')=(\cB\otimes\id_\eH)'=\cB'\otimes\sB(\eH)$ simply ``doubles'' the action of the ideal $\U_1\otimes\sB(\eH)$ and leaves the ideal $(\U-\U_1)\otimes\sB(\eH)$ fixed.
\eex

\setlength{\baselineskip}{2.5ex}



\newcommand{\Swap}[2]{#2#1}\newcommand{\Sort}[1]{}
\providecommand{\bysame}{\leavevmode\hbox to3em{\hrulefill}\thinspace}
\providecommand{\MR}{\relax\ifhmode\unskip\space\fi MR }
\providecommand{\MRhref}[2]{%
  \href{http://www.ams.org/mathscinet-getitem?mr=#1}{#2}
}
\providecommand{\href}[2]{#2}


\end{document}